\documentclass{article} % current default for manuscript submission

\usepackage[backref=true,backend=biber,giveninits=true,style=authoryear,uniquename=init,maxbibnames=99]{biblatex}

\usepackage[caption=false]{subfig}

\usepackage{mathtools}
\usepackage[ruled,vlined]{algorithm2e} % Algorithms
\usepackage{xspace}
\usepackage{amsmath}
\usepackage{bbm}

\usepackage{booktabs} % for professional tables
\usepackage{amssymb}
\usepackage{ifthen}
\usepackage{url}
\usepackage{graphicx}
\usepackage{color}
\usepackage{array}
\usepackage{theorem}

\usepackage[T1]{fontenc}
\usepackage[utf8]{inputenc}
\usepackage[english]{babel}
\usepackage{fullpage}

\usepackage[colorlinks]{hyperref}
\usepackage{hypernat}
\usepackage{xcolor}
\definecolor{NavyBlue}{RGB}{35,35,142}
\definecolor{RawSienna}{RGB}{199,97,20}
\hypersetup{
colorlinks,%
citecolor=NavyBlue,%
filecolor=NavyBlue,%
linkcolor=RawSienna,%
urlcolor=NavyBlue
}
\usepackage{authblk}

\usepackage[toc,page]{appendix}

\DeclarePairedDelimiter\floor{\lfloor}{\rfloor}

\newcommand{\bfx}{\mathbf{x}}
\newcommand{\bfy}{\mathbf{y}}
\newcommand{\bfz}{\mathbf{z}}

\newcommand{\bfI}{\mathbf{I}}

\renewcommand{\epsilon}{\varepsilon}

\newcommand{\calV}{\mathcal{V}}

\newcommand{\calA}{\mathcal{A}}

\newtheorem{proposition}{Proposition}[section]
\newtheorem{remark}{Remark}[section]
\newtheorem{result}{Result}[section]

\begin{document}
%%%%%%%%%%%%%%%%

\title{Solving a Continent-Scale Inventory Routing Problem at Renault}

\author[1,2]{Louis Bouvier\footnote{Corresponding author: \href{mailto:louis.bouvier@enpc.fr}{louis.bouvier@enpc.fr}}} 
\author[1]{Guillaume Dalle}
\author[1]{Axel Parmentier}
\author[3]{Thibaut Vidal}
\affil[1]{CERMICS, Ecole des Ponts, France}
\affil[2]{Groupe Renault}
\affil[3]{Polytechnique Montréal, Canada}

\setcounter{Maxaffil}{0}
\renewcommand\Affilfont{\itshape\small}

\maketitle

\begin{abstract}
  This paper is the fruit of a partnership with Renault.
  Their reverse logistic requires solving a continent-scale multi-attribute inventory routing problem (IRP).
  With an average of 30 commodities, 16 depots, and 600 customers spread across a continent, our instances are orders of magnitude larger than those in the literature. Existing algorithms do not scale, 
  so we propose a large neighborhood search (LNS). 
  To make it work, (1) we generalize existing split delivery vehicle routing problem and IRP neighborhoods to this context, (2) we turn a state-of-the-art matheuristic for medium-scale IRP into a large neighborhood, and (3) we introduce two novel perturbations: 
  the reinsertion of a customer and that of a commodity into the IRP solution. 
  We also derive a new lower bound based on a flow relaxation.
  In order to stimulate the research on large-scale IRP, we introduce a library of industrial instances.
  We benchmark our algorithms on these instances and make our code open-source. 
  Extensive numerical experiments highlight the relevance of each component of our LNS.
\end{abstract}

\

\section{Introduction}\label{intro} %%1.
% Commodities delivery 
% When a supplier manages the delivery of commodities to its customers on a multiple-day horizon in a centralized manner,
% the process is known as vendor-managed-inventory (VMI), opposed to  the retailer-managed-inventory (RMI) as detailed in \textcite{archettiInventoryRoutingProblem2016}.
% The optimization task to address then is an inventory routing problem (IRP).
The inventory routing problem (IRP) arises when a supplier manages the delivery of commodities to its customers on a multiple-day horizon in a centralized manner \parencite{archettiInventoryRoutingProblem2016}.
It consists in planning routes to deliver commodities from depots to customers with the objective of minimizing 
inventory and routing costs. 
% It is a NP-hard problem (it incorporates the split delivery vehicle routing problem (SDVRP) which is itself a NP-hard problem)
This NP-hard problem has received significant attention in the operations research literature over the past\:$40$ years.

The present paper is motivated by a partnership with Renault, a major European car manufacturer who must routinely solve IRP instances of unprecedented continental scale and complexity as part of their reverse logistic problem.
Indeed, they receive car parts from suppliers at their plants in packaging, and reuse the latter, which implies the need for reverse packaging logistics. 
The goal of our partnership is to redesign their IRP algorithm.
% involving IRP instances 
% of an unprecedented continent scale and complexity: 
This is challenging because of (1)\:the size of the resulting instances, with\:$600$ customers and\:$16$ depots on average, (2)\:the\:$30$ different commodities involved, and (3)\:the specific challenges that arise from the geography and timescale. 
When depots and customers are scattered across a whole continent, travel times can last up to ten days. Beyond requiring a long horizon, 21 days in our case, this makes the problem more difficult because classic decoupling results on the IRP, which were exploited in previous algorithms, are no longer valid. 
For instance, changing the order of the customers along a route impacts the arrival day at each customer and therefore the customer inventory levels. Hence, routes with suboptimal routing cost may be better because of inventory cost, which is not usually the case with the IRP. It can be compared with the continuous-time IRP discussed in \textcite{savelsberghOptimizationAlgorithmInventory2008} or \textcite{lagosContinuousTimeInventoryRoutingProblem2020}. \textcite{touzout_assign-and-route_2022} go further in this direction by considering a time-dependent setting in which travelling time depends on the departure time. We do not include this specificity in our work because congestion does not significantly impact long trips.
Finally, our partner's supply chain process requires that (4)\:the solution algorithm should not take more than\:$90$ minutes on our computing cluster.

% . The 
% multi-depot and multicommodity aspects entail an inherently hard problem.
State-of-the-art exact algorithms rely on branch-and-cut \parencite{archettiBranchandCutAlgorithmVendorManaged2007,coelhoBranchandcutAlgorithmMultiproduct2013,manousakisImprovedBranchandcutInventory2021} 
and branch-and-price-and-cut methods \parencite{desaulniersBranchPriceandCutAlgorithmInventoryRouting2015}
with dedicated valid inequalities.
They can optimally solve single-commodity single-depot instances with up to\:$50$ customers, but are not appropriate for our large-scale setting.

Typical heuristics include route-based matheuristics \parencite{fischettiMatheuristics2016}, decomposition matheuristics, and metaheuristics. In this field, \textcite{bertazziMatheuristicAlgorithmMultidepot2019} and \textcite{archettiMatheuristicMultivehicleInventory2017}
are route-based matheuristics.
The main idea is to reduce the size of the mixed-integer linear program (MILP) formulation of the IRP, by selecting promising routes heuristically. 
Although \textcite{bertazziMatheuristicAlgorithmMultidepot2019} is dedicated to the multi-depot case, neither of the two papers handles the multicommodity aspect we must face,
and their largest instances have up to six days horizon, six depots and\:$50$ customers. 
The methods cannot be applied directly in our context, because the MILP remains too large, even when we restrict ourselves to ``promising routes''. 
We instead adapt their principle to our setting, leading to the ``reload fixed-path vehicles'' subroutine. 

Another common approach is to tackle the IRP through a decomposition \parencite{campbellDecompositionApproachInventoryRouting2004,cordeauDecompositionbasedHeuristicMultipleproduct2015}. For instance, first set the quantities to be sent, and then create the routes to respect them.
The largest instances solved with this two-step method have a single depot, up to five commodities, a six days horizon and\:$50$ customers. 
We also adapt it to our setting and use it as an initialization heuristic.

Some metaheuristics have been designed for real-world IRP. For instance, \textcite{benoistRandomizedLocalSearch2011a} introduce a randomized local search to address a large-scale single-commodity IRP with pickups, time windows, driver safety and other constraints that are specific to their use case, but less relevant to ours. 
\textcite{suMatheuristicAlgorithmInventory2020} address a real-world IRP 
from the ROADEF IRP-Challenge\:$2016$ available at \url{https://www.roadef.org/challenge/2016/fr/} with a large neighborhood search 
based on mathematical programming. As for \textcite{benoistRandomizedLocalSearch2011a},
the single-commodity formulation with additional constraints is not adapted to our multicommodity context.
Other large neighborhoods are introduced in \textcite{nolzStochasticInventoryRouting2014} for 
a single-depot single-commodity stochastic version of the IRP. They are based on perturbation ideas such as the removal of every customer visit on a particular day, followed by a best-insertion policy. This process of removal and 
insertion is at the core of our work, but we leverage MILP formulations for the insertion. 
A kernel search heuristic based on a preliminary tabu search 
is considered in \textcite{archettiKernelSearchHeuristic2021}. In this framework, smaller MILPs with increasing size are solved iteratively to improve an initial solution. 
Single-depot, single-commodity instances with up to six days horizon and\:$200$ customers are solved. 
The study \textcite{coelhoVariableMIPNeighborhood2020} fixes a part of the decision variables, this time based on the problem's main ``axes'' -- that is to say the different types of sites 
involved in the deliveries, the routing and the inventory aspects. It solves reduced MILPs in a variable neighborhood search, and defines the multi-attribute IRP
as a multicommodity, multi-depot and multi-vehicle IRP. 
This constraint structure is the closest to ours. Nonetheless, the instances 
are smaller, with up to six days horizon,\:$50$ customers, six depots and three commodities.
Besides, this study is restricted to one-day routes, whereas we explicitly deal with routes that last multiple days, which creates an additional combinatorial challenge.
Another difference is that it incorporates a heterogeneous fleet with three distinct vehicle types, whereas we consider a homogeneous infinite fleet of vehicles. Other approaches follow this idea of fixing a part of the solution. For instance, \textcite{chitsazUnifiedDecompositionMatheuristic2019} use a three-phase decomposition matheuristic to address the assembly routing problem. The authors also derive ways to adapt it to the single-depot single-commodity IRP. Similarly, \textcite{vadsethIterativeMatheuristicInventory2021} introduce an iterative matheuristic for the single-depot single-commodity IRP. The method alternates between generating a small set of routes and solving a path-flow formulation given the fixed routes. Our approach shares this idea of updating sequentially the set of routes and the quantities delivered. Nonetheless, due to our multicommodity aspect, optimizing the quantities given the routes remains an $\mathcal{NP}$-hard problem in our case. We also introduce additional perturbations.
Therefore, to the best of our knowledge, no algorithm is known to properly scale to our instances. In this context, our main contributions are the following: 
\begin{enumerate}
  \item We introduce two new large-scale perturbation neighborhoods designed for the multi-attribute IRP. They are based on well-solved MILP formulations, and enable to escape from local minima.  
  \item We design an efficient large neighborhood search (LNS) built upon these large neighborhoods. We generalize\:$12$ Traveling Salesman Problem (TSP), and Split Delivery Vehicle Routing Problem (SDVRP) (see, e.g., \textcite{drorSplitDeliveryRouting1990a,archettiSplitDeliveryVehicle2008a}) neighborhoods from the literature to our IRP context.
  We propose a new matheuristic inspired by \textcite{archettiMatheuristicMultivehicleInventory2017} and \textcite{bertazziMatheuristicAlgorithmMultidepot2019} and adapted to our large-scale setting.
  \item We compute a new lower bound based on a linear program (LP) relaxation (one flow per commodity). To the best of our knowledge, this 
  relaxation is not considered in the literature. We do not expect the bound to be tight, which is a feature of every relaxation in the IRP literature.
  But it is useful to compare algorithm performance on instances with distinct scales.
  \item We provide a publicly available library of realistic multi-attribute IRP instances of a continent scale, as an incentive for further research on the topic. 
  \item We give access to our open-source Julia \parencite{Julia-2017} package that implements the ideas of the present paper. It is available on GitHub \parencite{bouvier_louis_2023_8179161}, and we also provide the library of instances \parencite{louis_bouvier_2023_8177237} and the solutions we obtained \parencite{louis_bouvier_2023_8177271}.
  \item We proceed to extensive numerical experiments. Since no algorithm is known to scale to our context, we compare the adapted route-based matheuristic and our large neighborhood search.
\end{enumerate}

We precisely define the problem we consider in Section\:\ref{problemdescription}. We then provide an overview 
of the different solution processes in Section\:\ref{sec:algos_overview}. The three next sections\:\ref{sec:flow_bin_pack}, \ref{sec:routing_LS} and\:\ref{sec:MILP_based} detail the algorithms. 
We finally report our numerical experiments in Section\:\ref{sec:numerical_results}.

\section{Problem Description}\label{problemdescription} %% 2.

%\subsection{General description}\label{subsec:general_description}
We consider a rich variant of the inventory routing problem.
A supplier manages the delivery of several commodities to its customers on a multiple-day horizon in a centralized manner \parencite{archettiInventoryRoutingProblem2016}.
The supplier has to plan routes to deliver commodities from its depots to customers with the objective of minimizing 
inventory and routing costs. Every day, depots release and customers demand commodities. Routes may last several days due to large travel distances. For each route, the starting depot, the ordered list of customers visited, as well as the quantities of each commodity to be delivered at each stop must be decided. We formalize this problem and include additional details in the remaining of Section\:\ref{problemdescription}. 

\begin{figure}[htb]
  \subfloat[Tiny IRP instance and corresponding solution.]{
  \begin{minipage}{0.7\textwidth}
    \centering
    \includegraphics[width=.99\linewidth]{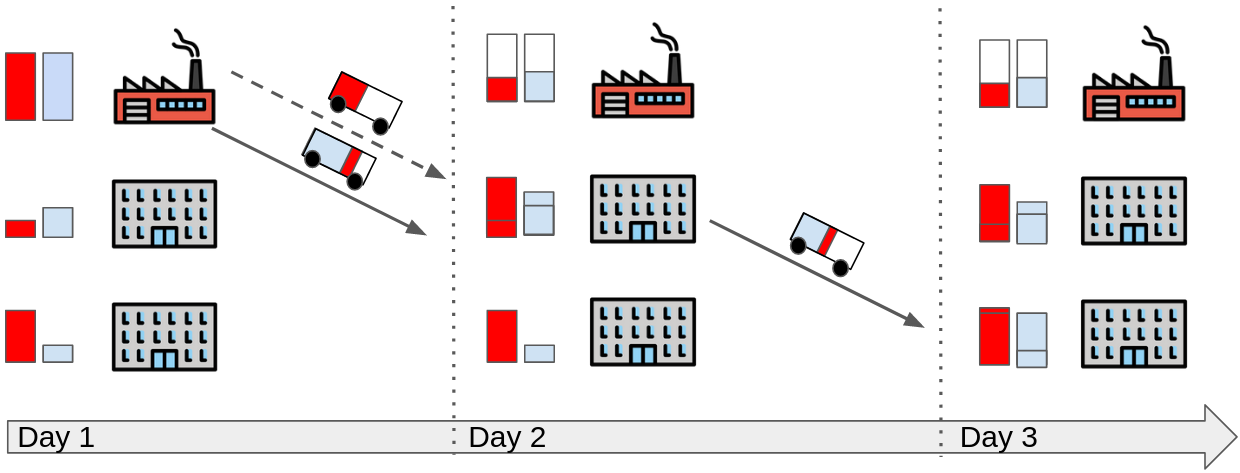}
  \end{minipage}
  }
  \subfloat[Real European solution.]{
  \begin{minipage}{0.3\textwidth}
    \centering
    \includegraphics[width=0.99\linewidth]{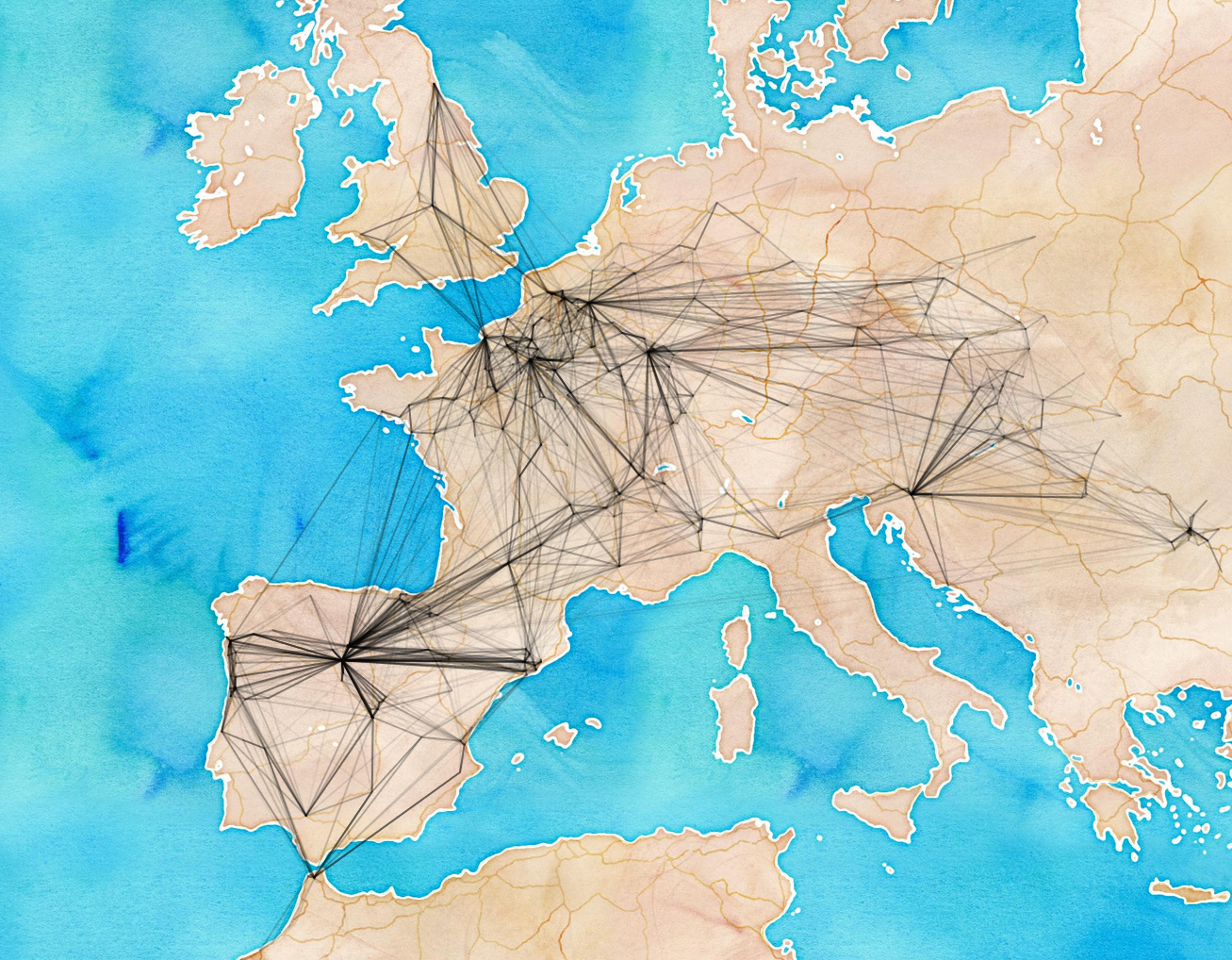}
  \end{minipage}
  }
  \caption{Instance and solution. Tiny example with two customers and two commodities (a), real European instance and routes solution in black (b).}
  \label{fig:IRP_illustration}
\end{figure}

Figure\:\ref{fig:IRP_illustration}\:(a) illustrates the routing, bin packing, and inventory aspects on a tiny example instance including one depot (in red with chimneys), two customers (in grey), two commodities (red and blue) and a three-days horizon. For simplicity, neither release nor demand occurs. Only initial inventories and deliveries are factored in the inventory dynamics.
Inventory levels are proportional to the sizes of the bars on the left of the sites, and measured in the evening of each day. The space occupied by the commodities loaded in the vehicles on each transport section is proportional to the bars above each arrow. 
On the first day, two routes start from the depot (one represented by the dotted arrow, the other by the full arrow) and reach the first customer on the second day. On the third day, the second route (full arrow) reaches the second customer. No route starts on the third day. Inventory levels are updated in accordance with the quantities delivered. To illustrate the complexity of our instances, Figure\:\ref{fig:IRP_illustration}\:(b) shows the solution of a real European instance provided by our algorithm. Routes are drawn with black lines on the continent. The thicker the line, the more the route is used. We understand well the continuous-time aspect at this scale.

\subsection{Notations and Data}
\label{subsec:notations}

Let\:$\mathbb{Z}^+$ be the set of non-negative integers.
For\:$a \in \mathbb{Z}^+$, we denote by\:$[a]$ the set\:$\{1,...,a\}$. 
Besides, for\:$x \in \mathbb{R}$, we define\:${(x)^+ := \max(x,0)}$. We denote by\:$|\mathcal{S}|$ the cardinal of a set or list\:$\mathcal{S}$.
When we explicitly consider a vector ${\bfx = (x_1,...,x_p)}$ of dimension ${p \in \mathbb{Z}^+}$, we use the notation\:${\bfx \in \mathbb{Z}}$ or\:${\bfx \in \{0,1\}}$ instead of\:${\bfx \in \mathbb{Z}^{|p|}}$ or\:${\bfx \in \{0,1\}^{|p|}}$ respectively.
Let\:$M$ be the set of \emph{commodities}, $D$ the set of \emph{depots} and $C$ the set of \emph{customers}, that respectively \emph{release} and \emph{demand} commodities\:$m \in M$. 
The \emph{time horizon} is\:$T \in \mathbb{Z}^+$ days.
At the beginning, each \emph{vertex}\:$v$ (depot or customer) has 
an \emph{initial inventory} of commodity\:$m$ denoted by\:$I^0_{mv}$. On each day\:$t \in [T]$, 
a customer\:$c$ demands a quantity\:$b^-_{mct}$ of commodity\:$m$. A depot\:$d$ releases 
a quantity\:$b^+_{mdt}$ of commodity\:$m$.
We say that a depot\:$d \in D$ \emph{uses} a commodity\:$m \in M$ if it has a positive initial 
inventory or a positive release for\:$m$ at least once over the horizon. We denote by\:$M_d$ the set 
of \emph{commodities used by depot}\:$d$. We similarly define\:$M_c$ as the set of \emph{commodities 
used by customer}\:$c \in C$, based on initial inventory and demand. 
A \emph{maximum inventory capacity}\:$\kappa_{mvt}$ is set on the night of each day\:$t$ per vertex\:$v$
and commodity\:$m$. Below this capacity, no inventory cost is paid. Above,
a cost is set to\:$c^{\texttt{exc}}_{mv}$ per unit, where ``\texttt{exc}'' stands for excess. Besides, a price\:$c^{\texttt{short}}_{mc}$ is paid 
per unit of unsatisfied demand for commodity\:$m$ of customer\:$c$, where ``\texttt{short}'' stands for shortage.
It corresponds to a soft constraint of non-negativity for the customers' inventories. 
We approximate commodities and vehicles by one-dimensional objects. 
We associate a length\:$\ell_m$ to each commodity\:$m \in M$.
We consider an infinite fleet of homogeneous vehicles of length\:$L$, to deliver the commodities from depots to customers.
They are not assigned to a particular depot. 
A 1D bin packing problem must be solved for vehicle loading. 
The depots and customers are the 
vertices\:$\mathcal{V} = D \cup C$ of a directed graph\:${\mathcal{D} = (\mathcal{V}, \mathcal{A})}$ that we name the \emph{locations graph}. 
The directed aspect is used to model the fact that transport durations and distances depend on the trip direction.
There is an arc\:$a = (u,v) \in \mathcal{A}$ for each vertex\:$u \in D \cup C$ and\:$v \in C$, $v \neq u$. Given a vertex $v$, we denote by\:$\delta^+(v)$ the set of arcs outgoing from $v$, and by\:$\delta^-(v)$ the set of arcs incoming to $v$. We associate 
a \emph{distance}\:$\Delta_a$ (in kilometers) and a \emph{transport duration}\:$\tau_a$ (in hours) to each arc.
We assume that the distances satisfy the triangular inequality.
When planning a route, a cost is paid per vehicle\:$c^{\texttt{veh}}$, per stop (customer visited)\:$c^{\texttt{stop}}$, and per 
kilometer travelled\:$c^{\texttt{km}}$. The number of stops must not exceed\:$S_{max}$, which is a practical requirement of the car manufacturer. The \emph{limit of driving hours 
per day} is\:$\tau_{\max}$. 
The IRP consists in building a set of routes (see Section\:\ref{subsec:routes}) 
to deliver commodities from depots to customers, minimizing the sum of the routing, inventory and shortage costs and respecting 
feasibility constraints detailed in Sections\:\ref{IRP_formulation} and\:\ref{subsec:routes}. 

\subsection{The Route Structure}
\label{subsec:routes}

An admissible path\:$P = (v_0, v_1, ..., v_k)$ in the locations graph is an elementary path, i.e., a path with pairwise distinct vertices. It starts from a depot\:$v_0 \in D$, and visits 
customers ${(v_1, ..., v_k) \in C^k, v_i \neq v_j}$. We emphasize it means routes are open, they do not end in their starting depots. We have a limit\:$S_{max}$ to the number of customers visited: 
\begin{equation}\label{eq:stop_number_limit}
  |P| \leq S_{max}+1.
\end{equation}

This constraint is required for practical reasons by the car manufacturer, because of the complexity of the unloading process. One alternative could be to limit the total duration or the total distance travelled by a route.
We highlight the fact that a path does not end at its starting depot. Let\:$\mathcal{P}$ be the set of admissible paths, and $A(P)$\:the set of arcs in a path\:$P$. A \emph{route}\:$r$ is a ``timed and loaded path''. It is a tuple\:${r = (t^r, P^r, \mathbf{q}^r)}$ where: 

\begin{itemize}
    \item $t^r \in [T]$ is the day of the departure. 
    \item $P^r = (v_0^r, v_1^r, ..., v_k^r) \in \mathcal{P}$ is the admissible path followed.
    \item $\mathbf{q}^r= (q^r_{ms})_{m \in M, s \in [|P^r|-1]} \in (\mathbb{Z}^+)^{|M|\times(|P^r|-1)}$ are the quantities delivered,
    for each commodity\:${m \in M}$ and to each customer\:$v^r_s$ for\:$s \in [|P^r|-1]$.
\end{itemize}

The total load must not exceed the vehicle capacity\:$L$, which can be written as: 
\begin{equation}\label{eq:durations_vehicle_capacity_constraint}
  \sum_{m \in M}\ell_m\bigg(\sum_{s \in [|P^r|-1]}q^r_{ms}\bigg) \leq L. 
\end{equation}

Given a route\:$r$, and the transport durations\:$\tau_a$ for\:$a \in \mathcal{A}$, we can compute the arrival day\:$t^r_s$ at the customer\:$v^r_s$ for\:$s \in [|P^r|-1]$ as follows.
We first compute the cumulated transport duration in hours up to customer\:$v_s^r$, with\:$\tau^r_0 = 0$: 
\begin{equation}\label{eq:route_duration}
  \tau_s^{r} = \tau^r_{s-1} + \tau_{(v_{s-1}^{r},v_{s}^{r})}, \quad \forall s \in [|P^r|-1].
\end{equation}
Then, the actual day\:$t_s^r$ of arrival at customer\:$v_s^r$ takes pauses into account: 
\begin{equation}\label{eq:date_arrival}
  t_s^r = t^r + \Big\lfloor\frac{\tau_s^{r}}{\tau_{\max}}\Big\rfloor, \quad \forall s \in [|P^r|-1].
\end{equation}
Equation \eqref{eq:date_arrival} means that when a driver exceeds the driving time limit per day\:$\tau_{\max}$, a pause is made until the next day. The vehicle then goes on from 
the location of the pause. Since in practice routes start from depots in the morning and the deliveries are only available in the evening at the customers, 
it is indeed a floor and not a ceiling function we consider in Equation \eqref{eq:date_arrival}.
A route must also visit every stop before the horizon\:$T$, which can be written as:
\begin{equation}\label{eq:dates_before_horizon}
  t^r_s \leq T, \quad \forall s \in [|P^r|-1].
\end{equation}

We henceforth denote by\:$\mathcal{R}$ the set of admissible routes.
A \emph{direct route} follows a path\:$P = (d,c)$ from a depot\:$d \in D$ to a customer\:$c \in C$ in the locations graph. It has only one 
arc. 

\subsection{Inventory Routing Formulation}\label{IRP_formulation}

The variables we consider are the following. We denote by\:$z_{mdt}^-$ the quantity of commodity\:$m$ sent from depot\:$d$ on day\:$t$, and by\:$z_{mct}^+$ the quantity of commodity\:$m$ delivered to customer\:$c$ on day\:$t$. Let\:$I_{mvt}$ be the inventory of commodity\:$m$ at vertex\:$v$ on the evening of day\:$t$.
We last denote by\:$x_r$ the number of vehicles following route\:$r$. We then consider the MILP formulation:
\begin{subequations}
  \begin{alignat}{2}
    \tag{multi-attribute-IRP}
  &\min_{\bfx, \bfz, \bfI} & \quad & \phantom{+}  \sum_{r}x_r \bigg(c^{\texttt{veh}} + c^{\texttt{stop}}(|P^r|-1) + c^{\texttt{km}} \sum_{a \in A(P^r)}\Delta_a \bigg)\\*
  &                        &       &  + \sum_{d,t,m} c^{\texttt{exc}}_{md} \left(I_{mdt} - \kappa_{mdt}\right)^+ \\*
  &                        &       &  + \sum_{c,t,m} c^{\texttt{exc}}_{mc}\left(I_{mct} - \kappa_{mct}\right)^+ + c^{\texttt{short}}_{mc} (b^-_{mct} -I_{mc(t-1)})^+\\*
  & \text{subject to} & & z_{mdt}^- = \sum_{\substack{r, \\v^r_0 = d, \\t^r = t}} \sum_{\substack{s \in [|P^r|-1]}}x_r q^r_{ms}, \quad \forall m \in M, \quad \forall d \in D, \quad \forall t \in [T] \label{eq:durations_q_to_z_depot}\\*
  & & &z_{mct}^+ = \sum_{r} \sum_{s \in [|P^r|-1], t^r_s = t, v^r_s = c}x_r q^r_{ms}, \quad \forall m \in M, \quad \forall c \in C, \quad \forall t \in [T] \label{eq:durations_q_to_z_cust}\\*
  & & &I_{mdt} =I_{md(t-1)}+b_{mdt}^+ -z_{mdt}^-, \quad \forall m \in M, \quad \forall d \in D, \quad \forall t \in [T] \label{eq:durations_inventory_dynamics_depots}\\*
  & & &I_{md0} =I_{md}^0, \quad \forall m \in M, \quad \forall d \in D \label{eq:durations_inventory_intial_depots}\\*
  & & &I_{mct} = \big(I_{mc(t-1)}-b_{mct}^-\big)^+ +z_{mct}^+, \quad \forall m \in M, \quad \forall c \in C, \quad \forall t \in [T] \label{eq:durations_inventory_dynamics_customer}\\*
  & & &I_{mc0} =I_{mc}^0, \quad \forall m \in M, \quad \forall c \in C \label{eq:durations_inventory_intial_customers}\\*
  % & & &I_{mct} = 0, \quad \forall c \in C, \quad \forall t \in \zeroT, \quad \forall m \notin M_c \label{eq:durations_control_cust_commodities}\\*
  & & &\bfx \geq 0, ~ \bfz \geq 0, ~ \bfI \geq 0 \label{eq:nonegativity}\\*
  & & &\bfx \in \mathbb{Z}, ~ \bfz \in \mathbb{Z}, ~ \bfI \in \mathbb{Z}\label{eq:integers}
\end{alignat}
\end{subequations}

We notice that, given an IRP instance and the route variables\:$\bfx$, we can deduce the quantities sent or received\:$\bfz$ and the inventory\:$\bfI$. The latter are useful to express the IRP as an 
MILP.

\paragraph{Objective function.}
The first sum models a cost per vehicle, per stop and per kilometer travelled. The second one is related to the excess inventory during the nights at the depots. The third one
has both an excess inventory and a shortage part. The quantity\:$\big(b_{mct}^- -I_{mc(t-1)}\big)^+$
is a substitute that is bought separately when a shortage appears.

\paragraph{Constraints.}
The\:$\bfx$ variable is used to count the number of vehicles that follow the admissible routes defined 
in Section\:\ref{subsec:routes}.
Equation \eqref{eq:durations_q_to_z_depot} is used to bind the total quantities that are sent from each depot to the route deliveries at each customer.
Equation \eqref{eq:durations_q_to_z_cust} links the total quantities received per customer to the route deliveries. 
Constraints \eqref{eq:durations_inventory_dynamics_depots}-\eqref{eq:durations_inventory_intial_depots} define the inventory dynamics at the depots, and
\eqref{eq:durations_inventory_dynamics_customer}-\eqref{eq:durations_inventory_intial_customers} at the customers. 
We highlight we cannot deliver a commodity to a customer that does not need it -- in the sense of $C_m$ -- because the maximum inventory capacity is set to zero and the excess inventory cost to infinity. The MILP (multi-attribute-IRP) is intractable over our instances, we instead suggest several heuristic approaches in the next section.

\section{Overview of the Algorithms and General Concepts}\label{sec:algos_overview}
We emphasize the main principles of the algorithms in Section\:\ref{subsec:algos} before going into the details of each of their components. We also introduce generic flow graphs and formulations in Section\:\ref{sec:flow_graphs_formulations}, concepts useful in the rest of the present paper.

\subsection{Overview of the Algorithms}\label{subsec:algos}
We compare three algorithms with increasing degrees of sophistication and performance: an \emph{initialization + local search} algorithm to quickly derive non-trivial IRP solutions, 
a \emph{route-based matheuristic}, and our LNS. The two first algorithms are adapted from frameworks of the literature, the last one is our main contribution.

\subsubsection{Subroutines}\label{subsec:subroutines} 

\begin{figure}[htb]
  \centering
  \includegraphics[width=.99\linewidth]{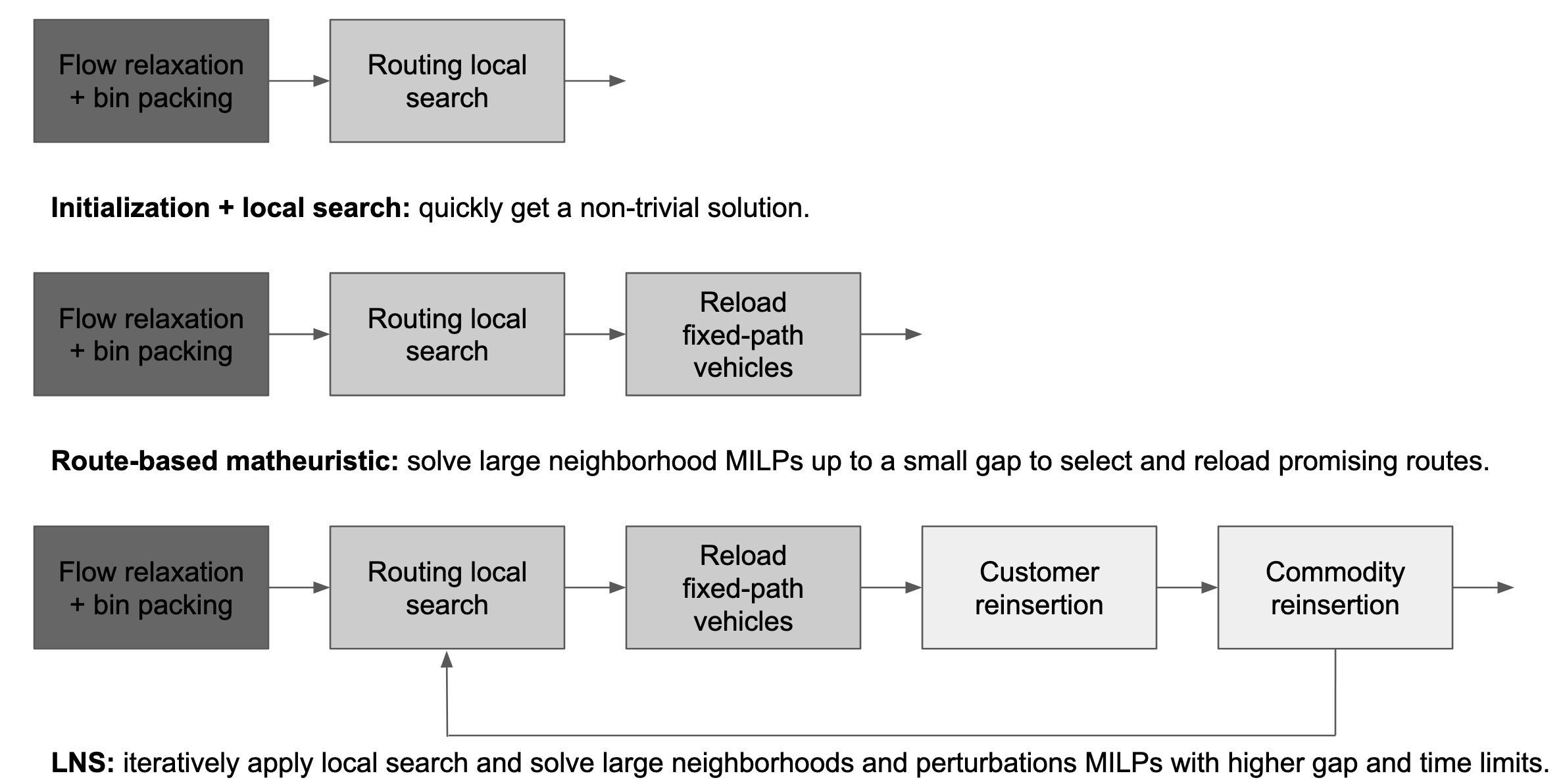}
  \caption{Different algorithms ordered by degree of sophistication.
  Dark gray corresponds to initialization, gray to descent subroutines, and light gray to perturbation subroutines.}
  \label{fig:algorithms}
\end{figure}

Our algorithms are illustrated on Figure\:\ref{fig:algorithms} and combine five subroutines. We call inner iteration an iteration within any subroutine, and outer iteration a path through the four types of neighborhoods in the LNS (see the loop in Figure\:\ref{fig:algorithms}). The first subroutine builds an initial solution.  1)\:The \emph{flow relaxation + bin packing} (\texttt{flow relaxation}, \texttt{bin packing}) subroutine solves a flow relaxation -- thus an LP -- per commodity and deduces direct routes by approximately solving bin packing problems to respect vehicle capacity. The other four subroutines improve or perturb an existing solution, and can be applied any number of times in any order. 2)\:The \emph{routing local search} ($\texttt{routing local search}$) subroutine takes a random subset of routes and applies a local search with TSP and SDVRP neighborhoods (see Section\:\ref{subsec:route_moves}). 3)\:The \emph{reload fixed-path vehicles} ($\texttt{reload fixed path vehicles}$) subroutine solves an MILP per depot to re-optimize the load of the routes starting from it. The MILP is solved with a very low gap threshold in the route-based matheuristic, and up to a larger gap threshold in the LNS, adding a time limit. 4)\:The \emph{customer reinsertion} ($\texttt{customer reinsertion}$) subroutine removes a customer from every delivery of a solution and solves an MILP to reinsert it in the existing routes, also creating new direct routes. 5)\:The \emph{commodity reinsertion} ($\texttt{commodity reinsertion}$) subroutine removes a commodity from every delivery of a solution and solves an MILP to reinsert it.
These MILPs are solved up to a gap and time limit. We see the last four subroutines as local search procedures, and call them with a customizable number of inner iterations per outer iteration. In Section\:\ref{sec:flow_bin_pack}, we detail the flow relaxation + bin packing subroutine. The routing local search subroutine is described in Section\:\ref{sec:routing_LS}. The three large neighborhoods -- reload fixed-path vehicles, customer and commodity reinsertion -- are detailed in Sections\:\ref{reload}, \ref{customer_reinsertion} and\:\ref{commodity_reinsertion} respectively.

\subsubsection{Algorithms}\label{subsec:algorithms} 

\paragraph{Initialization + local search.}
This algorithm simply runs the initialization + bin packing subroutine to build an initial solution, and then applies the routing local search subroutine to improve it. It has the advantage of being fast (about four minutes on our large-scale instances on average) since it is based on an LP. It is detailed in Algorithm\:\ref{algo:greedy_heuristic}.

\begin{algorithm}[!h]
  \SetKwData{Action}{action}
  \SetKwFunction{multipleflowlp}{$\texttt{flow relaxation}$}\SetKwFunction{firstfitdecreasing}{$\texttt{bin packing}$}\SetKwFunction{multidepotls}{$\texttt{routing local search}$}
  \SetKwInOut{Input}{input}\SetKwInOut{Output}{output}
  \Input{$\mathcal{I}$ an IRP instance.}
  \Output{A solution $\textbf{r}$ to the IRP.}
  $\bfy$ = \multipleflowlp($\mathcal{I}$)\;
  $\textbf{r}$ = \firstfitdecreasing($\mathcal{I}$, $\bfy$)\;
  $\textbf{r}$ = \multidepotls($\mathcal{I}$, $\textbf{r}$)\;
  \caption{Initialization + local search}\label{algo:greedy_heuristic}
\end{algorithm}

\paragraph{Route-based matheuristic.}
As discussed in Section\:\ref{intro}, \textcite{archettiMatheuristicMultivehicleInventory2017, bertazziMatheuristicAlgorithmMultidepot2019}
solve the IRP given a subset of ``promising routes'' that are defined heuristically. This corresponds to reducing the set of feasible solutions to the IRP, 
which allows solving an MILP. The routes are either selected among those created during 
a tabu search and leading to cost improvements, or in a constructive manner. The underlying assumption is that they are likely to appear in a good IRP solution.
We adapt this idea to our 
setting. The main difficulty is that our MILP (multi-attribute-IRP) is intractable, even
when we restrict the set of ``promising routes'' to the set of an initial solution. 
This is due to the multicommodity aspect, the scale of our instances and the routes that last several days. We solve the restricted MILP heuristically, using a large neighborhood approach: 1)\:Apply the initialization + local search algorithm to get an initial solution.
2)\:Take the current solution as set of promising routes, and solve sequentially one MILP per depot, with the corresponding promising routes
that start from it. It is detailed in Algorithm\:\ref{algo:route_based_matheuristic}.

\begin{algorithm}[!h]
  \SetKwFunction{commodityreinsertionMILP}{$\texttt{commodity insertion MILP}$}\SetKwFunction{fillformerroutes}{$\texttt{fill former routes}$}\SetKwFunction{fillnewroutes}{$\texttt{fill new routes}$}\SetKwFunction{singledepotLS}{$\texttt{single depot local search}$}
  \SetKwInOut{Input}{input}\SetKwInOut{Output}{output}
  \Input{$\mathcal{I}$ an IRP instance, $\textbf{r} = (r_k)_{1 \leq k \leq K}$ the current solution with $K \in \mathbb{Z}^+$ routes, $\texttt{time limit}$ a time limit.}
  \Output{The solution $\textbf{r}$ updated.}
    $\textbf{r} = \texttt{initialization + local search}(\mathcal{I})$\;
    \For{$d \in D$}{
      $\textbf{r} = \texttt{reload fixed path vehicles}(\mathcal{I}, \textbf{r}, d)$; \tcp*[f]{small gap}\

      \If{\texttt{time elapsed} $\geq \texttt{time limit}$}
      {break;}
    }
  \caption{Route-based matheuristic}\label{algo:route_based_matheuristic}
\end{algorithm}

\paragraph{Large neighborhood search.}
The large neighborhood search Algorithm\:\ref{algo:LNS} first uses the initialization + local search approach to find a good initial solution. It then explores four kinds of neighborhoods.
Two of them always improve the solution: the routing local search and reload fixed-path vehicles subroutines. Contrary to the route-based matheuristic, the latter is applied with a greater gap limit and an additional time limit, in order to avoid 
spending too much time within it and cycle over the neighborhoods instead.
The two remaining ones are perturbations, which means they can deteriorate the solution.
They fix a part of the current solution and optimize the quantities and routes involving a particular customer 
or commodity over the entire horizon. They both lead to substantial changes, allowing the search to escape from 
local minima. The LNS uses both iteratively, selecting every customer and ten commodities at random per outer step. Those hyperparameters are tuned experimentally as shown in Appendix\:\ref{app:hyperparameters}.
The outer LNS iterations are illustrated by the loop arc on Figure\:\ref{fig:algorithms}. The LNS returns the best solution found, comparing after each subroutine the current solution with the best one so far.
The main idea behind this LNS is to consider the structure of the IRP, ``decompose'' it along its major axes, and solve smaller 
natural problems to explore the solution space. One difference with \textcite{coelhoVariableMIPNeighborhood2020} or 
\textcite{archettiKernelSearchHeuristic2021}, is that our idea is not to fix a part of the MILP variables (multi-attribute-IRP) and optimize with respect to the remaining ones, but to define new smaller MILPs based on the structure of the IRP. Let us now introduce a concept that helps describing our subroutines.

\begin{algorithm}[!h]
\SetKwFunction{commodityreinsertionMILP}{$\texttt{commodity insertion MILP}$}\SetKwFunction{fillformerroutes}{$\texttt{fill former routes}$}\SetKwFunction{fillnewroutes}{$\texttt{fill new routes}$}\SetKwFunction{singledepotLS}{$\texttt{single depot local search}$}
\SetKwInOut{Input}{input}\SetKwInOut{Output}{output}
\Input{$\mathcal{I}$ an IRP instance, $\textbf{r} = (r_k)_{1 \leq k \leq K}$ the current solution with $K \in \mathbb{Z}^+$ routes, $\texttt{time limit}$ a time limit, $n_{\texttt{cust}}$ the number of customers to reinsert, $n_{\texttt{comm}}$ the number of commodities to reinsert.}
\Output{The solution $\textbf{r}$ updated.}
$\textbf{r} = \texttt{initialization + local search}(\mathcal{I})$\;
\While{\texttt{time elapsed} $< \texttt{time limit}$}{
  $\textbf{r} = \texttt{routing local search}(\mathcal{I}, \textbf{r})$\;
  \For{$d \in D$}{
    $\textbf{r} = \texttt{reload fixed path vehicles}(\mathcal{I}, \textbf{r}, d)$; \tcp*[f]{large gap}\\
  }
  $C_{\texttt{sub}} = sample(C, n_{\texttt{cust}})$\;
  \For{$c \in C_{\texttt{sub}}$}{
    $\textbf{r} = \texttt{customer reinsertion}(\mathcal{I}, \textbf{r}, c)$; \tcp*[f]{large gap}\\
  }
  $M_{\texttt{sub}} = sample(M, n_{\texttt{comm}})$\;
  \For{$m \in M_{\texttt{sub}}$}{
    $\textbf{r} = \texttt{commodity reinsertion}(\mathcal{I}, \textbf{r}, m)$; \tcp*[f]{large gap}\
  }
}
\caption{Large neighborhood search}\label{algo:LNS}
\end{algorithm}

\subsection{Flow Graphs and Formulations}\label{sec:flow_graphs_formulations}

Let us consider the generic MILP formulation:
\begin{subequations}
  \begin{alignat}{2}\tag{generic-flow-MILP}
    &\min_{\bfy, \bfx}  & \quad & \phantom{+} \sum_{m \in M} \bfy_m^{\top} \textbf{c}_m + \sum_{r \in \mathbf{r}}x_rc_r\\*
    & \text{subject to} &       & \sum_{a \in \delta^+(v)} y_{ma} = \sum_{a \in \delta^-(v)} y_{ma}, \quad \forall m \in M, \quad \forall v \in \calV^m\label{eq:cirucaltion_commodities}\\*
    &                   &       & \bfy_m^{\min} \leq \bfy_m \leq \bfy_m^{\max}, \quad \forall m \in M \label{eq:capacity_commodity} \\*
    &                   &       &\sum_{m \in M} y_{ma}\ell_m \leq x_rL, \quad \forall a = (d \to (c, r)), \quad \forall r \in \mathbf{r} \label{eq:capa_routes}\\*
    &                   &       &\bfy \in \mathbb{Z} \label{eq:integer_flows}\\*
    &                   &       &\bfx \in \{0, 1\}. \label{eq:indicator_vehicle}
  \end{alignat}
  \end{subequations}

\begin{figure}[htb]
  \subfloat[Global structure]{
  \begin{minipage}{\textwidth}
    \centering
    \includegraphics[width=0.6\linewidth]{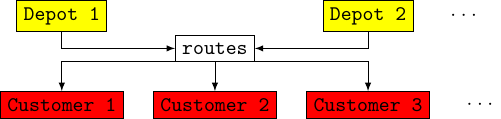}
  \end{minipage}
  }
  \newline
  \subfloat[Depot subgraph]{
  \begin{minipage}{\textwidth}
    \vspace{0.5cm}
    \centering
    \includegraphics[width=0.95\linewidth]{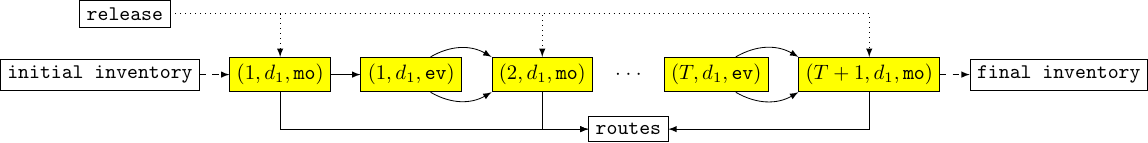}
  \end{minipage}
  }
  \newline
  \subfloat[Customer subgraph]{
     \begin{minipage}{\textwidth}
      \vspace{0.5cm}
    \centering
    \includegraphics[width=0.95\linewidth]{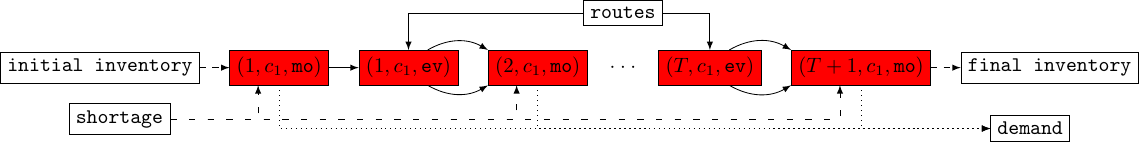}
  \end{minipage}
  }
  \caption{Commodity flow graph. Overview of the global graph structure (a), and then details of the subgraphs ``$\texttt{Depot 1}$'' (b) and ``$\texttt{Customer 1}$'' (c). The abbreviations ``$\texttt{mo}$'' and ``$\texttt{ev}$'' stand for morning and evening respectively.}
  \label{fig:overview_graph}
\end{figure}

The variable\:$\bfy_m$ encodes a flow on a given \emph{commodity graph} thanks to Equations \eqref{eq:cirucaltion_commodities}-\eqref{eq:capacity_commodity}. This flow enables modelling the depot and customer inventory dynamics of commodity\:$m$ defined by constraints\:\eqref{eq:durations_inventory_dynamics_depots}-\eqref{eq:durations_inventory_intial_customers}, and the quantities sent and received. Variable\:$\bfx$ encodes the routes that are used to deliver the commodities. Constraint \eqref{eq:capa_routes} indeed enforces the flows to respect the vehicle capacity when a route is used. The generic notation\:$(d \to (c, r))$ refers to the transport arc from the depot $d$ to the customer $c$ in the route $r$, on its departure day. In order to obtain a specific MILP formulation from this generic one, we must specify which commodity graph is used, and which set of routes\:$\mathbf{r}$ is considered.
We are going to use several distinct commodity graphs, see for instance the figures\:\ref{fig:commodity_graph}, \ref{fig:commodity_graph_refill}, \ref{fig:commodity_graph_customer_reinsertion} for detailed $(d \to (c, r))$ arcs. However, they all share a common structure which we describe now:

\paragraph{Depot subgraphs.} A subgraph per depot\:$d$ (Figure\:\ref{fig:overview_graph} (b)), which models its inventory dynamics. It is shared by the distinct formulations we introduce in this paper and has the following vertices: $(t,d,\texttt{morning})$ for\:$t \in [T+1]$, and\:$(t,d,\texttt{evening})$ for\:$t \in [T]$.

\paragraph{Customer subgraphs.} A subgraph per customer\:$c$ (Figure\:\ref{fig:overview_graph} (c)), which models its inventory dynamics. It is also shared by the distinct formulations we introduce in this paper, and it has the following vertices: $(t,c,\texttt{morning})$ for\:$t \in [T+1]$, and\:$(t,c,\texttt{evening})$ for\:$t \in [T]$.

\paragraph{A route subgraph.} Its specific structure depends on the formulation we consider. It contains paths between vertices of the form\:$(t,d,\texttt{morning})$ and vertices of the form\:$(\tilde{t},c,\texttt{evening})$ as shown on Figure\:\ref{fig:overview_graph} (a). Flow variables\:$\bfy$ on those paths model the quantities sent from depots to customers.

\paragraph{Artificial vertices.} In order to model commodity flows as circulations over commodity graphs, we add artificial vertices connected to the subgraphs above: 
\:$\texttt{source}$, $\texttt{sink}$, $\texttt{initial inventory}$, $\texttt{final inventory}$, $\texttt{release}$, $\texttt{shortage}$, and $\texttt{demand}$.

\paragraph{} The details of the arcs of the shared subgraphs defined above are in Table\:\ref{tab:arcs_generic_subgraphs}. For each arc, we give the following information. The subgraph it belongs to is first given. Then, we distinguish ``incoming'', ``outgoing'' and ``internal'' arcs with respect to the depots and customer subgraphs. A short description is stated to understand the meaning of the arcs. We specify the origin and destination vertices, as well as the minimum and maximum flow capacities associated to the flow variables on the arcs. Last, the cost corresponding to these variables are also given.

\paragraph{}
Given a graph\:$\mathcal{\tilde{D}} = (\tilde{\calV}, \tilde{\calA})$ with capacities associated to its arcs\:$(\bfy^{\min}, \bfy^{\max})$, we define the \emph{set of circulations} as
${\mathcal{C}(\mathcal{\tilde{D}}, \bfy^{\min}, \bfy^{\max}) = \{\bfy \in \mathbb{R}^{\tilde{\calA}}, \bfy_{\min} \leq \bfy \leq \bfy_{\max}, \forall v \in \tilde{\calV}, \sum_{a \in \delta^+(v)} \bfy_{a} = \sum_{a \in \delta^-(v)} \bfy_{a}\}}$. We use this notation instead of constraints \eqref{eq:cirucaltion_commodities}-\eqref{eq:capacity_commodity} in the rest of the paper. The minimum\:$\bfy^{\min}$ and maximum\:$\bfy^{\max}$ capacities are defined by the \textbf{Min} and \textbf{Max} columns in tables\:\ref{tab:arcs_generic_subgraphs}, \ref{tab:additional_arcs_commodity_graph_commodity_reinsertion} and figures\:\ref{fig:commodity_graph}\:(b), \ref{fig:commodity_graph_refill} (b) and \ref{fig:commodity_graph_customer_reinsertion} (b).

\begin{table}
  \caption{Arcs of the commodity flow graph shared by our formulations. \\ When not stated in the table, \textbf{Min} is\:$0$, \textbf{Cost} is\:$0$ and \textbf{Max} is\:$\infty$.}
  \setlength{\tabcolsep} {0.3cm} 
  \centering
\scalebox{0.65} 
{ 
\begin{tabular}{cccccccc} 
\toprule
\textbf{Subgraph} & \textbf{Arc type} &
\textbf{Arc description}     & \textbf{Origin}   & \textbf{Destination} & \textbf{Min} & \textbf{Max} & \textbf{Cost} \\
\midrule
Depot & Incoming & Initial inventory depot        & $\texttt{initial inventory}$       & $(1, d, \texttt{morning})$     & $I^0_{md}$  & $I^0_{md}$  &                                                       \\ 
Depot & Outgoing  & Final inventory depot          & $(T+1, d, \texttt{morning})$         & $\texttt{final inventory}$          &              &              &                                                       \\ 
Depot & Incoming  & Release depot                  & $\texttt{release}$              & $(t, d, \texttt{morning})$     & $b^+_{mdt}$  & $b^+_{mdt}$  &                                                       \\ 
Customer & Incoming  & Initial inventory customer     & $\texttt{initial inventory}$       & $(1, c, \texttt{morning})$     & $I^0_{mc}$  & $I^0_{mc}$  &                                                       \\ 
Customer & Outgoing  & Final inventory customer       & $(T+1, c, \texttt{morning})$         & $\texttt{final inventory}$          &              &              &                                                       \\ 
Customer & Outgoing  & Demand customer                & $(t, c, \texttt{morning})$         & $\texttt{demand}$           & $b^-_{mct}$  & $b^-_{mct}$  &                                                       \\  
Customer & Incoming  & Shortage customer              & $\texttt{shortage}$         & $(t, c, \texttt{morning})$     &              &              & $c^{\texttt{short}}_{mc}$                               \\ 
\midrule
Routes & & Transport $d \to c$            &    Formulation specific   &  (see Sections \ref{min_cost_relaxation}, \ref{reload}, \ref{customer_reinsertion}, \ref{commodity_reinsertion})      &              &                &   \\ 
\midrule
Depot  & Internal & Daily inventory depot          & $(t, d, \texttt{morning})$         & $(t, d, \texttt{evening})$      &              &              & \\ 
Depot  & Internal & Free night inventory depot     & $(t, d, \texttt{evening})$         & $(t+1, d, \texttt{morning})$   &              & $\kappa_{mdt}$    &                                                       \\ 
Depot  & Internal & Excess night inventory depot   & $(t, d, \texttt{evening})$         & $(t+1, d, \texttt{morning})$   &              &              & $c^{\texttt{exc}}_{md}$                                 \\ 
Customer  & Internal & Daily inventory customer       & $(t, c, \texttt{morning})$         & $(t, c, \texttt{evening})$      &              &              &   \\ 
Customer  & Internal & Free night inventory customer  & $(t, c, \texttt{evening})$         & $(t+1, c, \texttt{morning})$   &              & $\kappa_{mct}$    &                                                       \\ 
Customer  & Internal & Excess night inventory customer  & $(t, c, \texttt{evening})$       & $(t+1, c, \texttt{morning})$   &              &              & $c^{\texttt{exc}}_{mc}$                                 \\ 
\midrule
Artificial  &   & Circulation               & $\texttt{source}$       & $\texttt{release}$   &              &              &                                  \\ 
Artificial  &   & Circulation               & $\texttt{source}$       & $\texttt{initial inventory}$   &              &              &                                  \\ 
Artificial  &   & Circulation               & $\texttt{source}$       & $\texttt{shortage}$   &              &              &                                  \\ 
Artificial  &   & Circulation               & $\texttt{demand}$       & $\texttt{sink}$   &              &              &                                  \\ 
Artificial  &   & Circulation               & $\texttt{final inventory}$       & $\texttt{sink}$   &              &              &                                  \\ 
Artificial  &   & Circulation               & $\texttt{sink}$       & $\texttt{source}$   &              &              &                                  \\ 
\bottomrule
\end{tabular} %
}\label{tab:arcs_generic_subgraphs}
\end{table}

The (generic-flow-MILP) formulation is not fully specified here as we still need to define the route subgraph to obtain a formulation that can be given to a solver. It is an exact formulation if we enumerate all possible routes, but such a column generation is intractable on our large instances. Therefore, we derive heuristic formulations with specific route subgraphs in Sections\:\ref{min_cost_relaxation}, \ref{reload}, \ref{customer_reinsertion}, and \ref{commodity_reinsertion}.

\begin{remark}
  In (generic-flow-MILP), routes are modelled with individual paths in the commodity graphs, bound with indicator variables\:$\bfx$. Sometimes, we define them as paths in another flow graph. In this case, some vertices and arcs are shared between routes. We detail this aspect in Sections\:\ref{customer_reinsertion} and\:\ref{commodity_reinsertion}.
\end{remark}

\begin{remark}
  We sparsify the commodity graph. Instead of considering the sets of depots\:$D$ and customers\:$C$ in the commodity graph\:$\mathcal{D}^m$, we define the subsets\:${D_m = \{d \in D, m \in M_d\}}$ and\:${C_m = \{c \in C, m \in M_c\}}$ the depots and customers that use commodity\:$m$. We then restrict the depots and customer subgraphs of\:$\mathcal{D}^m$ to the ones of\:$D_m$ and\:$C_m$.
\end{remark}

\section{Flow Relaxation + Bin Packing Subroutine}\label{sec:flow_bin_pack}

The flow relaxation + bin packing subroutine is a fast heuristic to get an initial solution to (multi-attribute-IRP). It takes as input an IRP instance, and returns an initial IRP solution built from intermediate flow solutions that encode who sends what to whom and when.

\subsection{Multiple Minimum Cost Flows and Relaxation}\label{min_cost_relaxation}

\begin{figure}[htb]
  \centering
  \subfloat[Graph details]{
\begin{minipage}{1.\textwidth}
  \includegraphics[width=\linewidth]{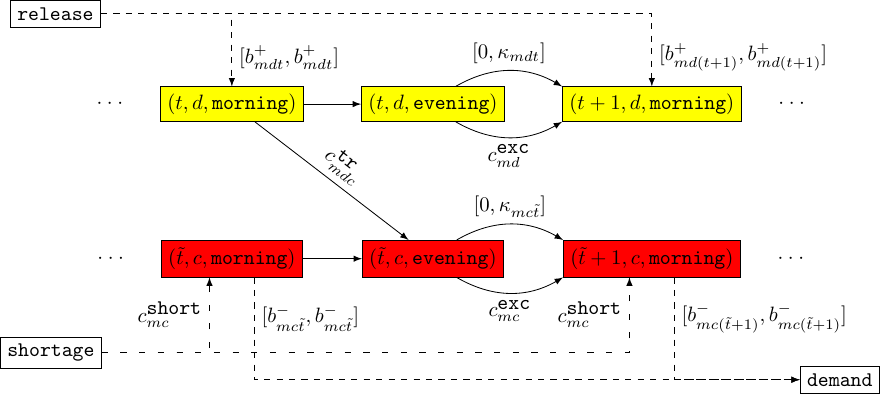}
\end{minipage}
  }
  \newline
  \subfloat[Additional commodity graph arcs compared with Table\:\ref{tab:arcs_generic_subgraphs}]{
\begin{minipage}{1.\textwidth}
  \vspace{0.5cm}
  \centering
  \setlength{\tabcolsep} {0.3cm} 
  \scalebox{0.8} 
  { 
  \begin{tabular}{ccccccc} 
  \toprule
  \textbf{Subgraph}  & \textbf{Arc description}     & \textbf{Origin}   & \textbf{Destination} & \textbf{Min} & \textbf{Max} & \textbf{Cost} \\
  \midrule
  Routes  & Transport $d \to c$            & $(t, d, \texttt{morning})$         & $(t+\floor{\frac{\tau_{dc}}{\tau_{\max}}}, c, \texttt{evening})$      &              &              &  $c^{\texttt{tr}}_{mdc}$\\ 
  Routes  & Delayed transport $d \to c$    & $(t, d, \texttt{morning})$         & $(\tilde{t}, c, \texttt{evening})$      &              &              &  $c^{\texttt{tr}}_{mdc}$\\ 
  \bottomrule
  \end{tabular} %
  }
\end{minipage}
  }
\caption{Details of the commodity graph\:$\mathcal{D}^m$ for the flow relaxation problem.}
\label{fig:commodity_graph}
\end{figure}

\paragraph{Minimum cost flow formulation.} Let\:$\bfy = (y_{ma})_{m \in M, a \in \calA^m}$ be a flow variable 
and\:$\textbf{c}$ be the corresponding costs defined
in Table\:\ref{tab:arcs_generic_subgraphs} for the shared subgraphs and Figure\:\ref{fig:commodity_graph} for the specific route subgraph. We consider the following LP:
\begin{subequations}
\begin{alignat}{2}
  \tag{flow-relaxation}
    &\min_{\bfy} {} & \quad & \phantom{+} \sum_{m \in M} \bfy_m^{\top} \textbf{c}_m  \\*
    & \text{subject to }&   & \bfy_m \in \mathcal{C}(\mathcal{D}^m, \bfy_m^{\min}, \bfy_m^{\max}), \quad \forall m \in M%\sum_{a \in \delta^+(v)} \bfy_{ma} = \sum_{a \in \delta^-(v)} \bfy_{ma}, \quad \forall m \in M, \quad \forall v \in \calV^m \\*
    %&                   &   & \bfy_m^{\min} \leq \bfy_m \leq \bfy_m^{\max} \quad \forall m \in M
\end{alignat}
\end{subequations}

In this variant of the generic MILP (Section\:\ref{sec:flow_graphs_formulations}), we do not introduce a route variable\:$\bfx$ and the corresponding cost. We instead consider a soft version of constraint\:\eqref{eq:capa_routes} in the commodity cost\:$\textbf{c}_m$. This leads to a separate flow LP for each commodity in (flow-relaxation). Let us now introduce the details of the commodity graph.

\paragraph{Details of the commodity graph.}
We define one graph per commodity\:$m$, named\:${\mathcal{D}^m = (\calV^m, \calA^m)}$. The vertices\:$\calV^m$ are exactly the ones defined in Section\:\ref{sec:flow_graphs_formulations}. The arcs detailed in Table\:\ref{tab:arcs_generic_subgraphs} are included. We specify the route subgraph in the table of Figure\:\ref{fig:commodity_graph}. It is made of direct transport and delayed transport\:$d \to c$ arcs. These arcs are added when the date of arrival at the customer is smaller than the horizon\:$T$. For each tuple\:$(t,d,c) \in [T] \times D_m \times C_m$, we add one delayed arc\:$\big((t, d) \to (\tilde{t}, c)\big)$ 
per possible delayed arrival day\:$\tilde{t}$ induced by an indirect path from\:$d$ to\:$c$ (thus visiting any set of other customers before\:$c$) that respects the route constraints defined by\:$\mathcal{R}$. 
Those possible delays are pre-computed, using a breadth-first search algorithm over the locations graph, with maximum depth set to\:$S_{\max}$. Indeed, we can browse the locations graph starting from depots, saving the cumulative delay at any vertex and any depth smaller than\:$S_{\max}$.
The flow on these arcs models the quantity of commodity sent from depot\:$d$ on day\:$t$ to customer\:$c$ with arrival on day\:$\tilde{t}$.
The precise structure of the graph is illustrated on Figure\:\ref{fig:commodity_graph}.
In arc annotations, capacities are given between brackets (e.g. $[0, \kappa_{mdt}]$), and costs without (e.g. $c_{md}^{\texttt{exc}}$). Dotted arrows 
are related to the shared artificial vertices, continuous ones to depots, customers and route subgraphs. When not stated in the table, \textbf{Min} is\:$0$, \textbf{Cost} is\:$0$ and \textbf{Max} is\:$\infty$. On this figure, only two days, one depot and one customer are shown. Besides, some artificial vertices are omitted for simplicity. The cost\:$c^{\texttt{tr}}_{mdc}$ is detailed below. 

%%%%
Since the routing price is paid at the vehicle level, we cannot derive a 
minimum cost commodity flow that takes it into account exactly without adding variables for each individual vehicle. Instead, 
we want to approximate this cost with transportation 
arcs between depots and customers
naturally involving commodity flow variables. 
A way to do so is to use a ``vehicle fraction'' unit per commodity, leading 
to:

\begin{equation}\label{eq:fraction_vehicle}
    c^{\texttt{tr}}_{mdc} = \frac{\ell_m}{L} (c^{\texttt{veh}} + c^{\texttt{stop}} + c^{\texttt{km}} \Delta_{dc}),\quad \forall d \in D_m,\quad \forall c \in C_m.
\end{equation}
 
In Equation \eqref{eq:fraction_vehicle} the factor\:$\frac{\ell_m}{L}$ is a way to scale the price paid for the delivery of a unit of commodity\:$m$ based on the percentage of 
a vehicle it occupies, hence the ``vehicle fraction''.

\begin{proposition} The optimization problem (flow-relaxation) based on\:$|M|$ flows is a relaxation 
of (multi-attribute-IRP). The optimal value of (flow-relaxation) is a lower bound to the cost of an optimal solution to our initial problem.
\end{proposition}

We now sketch the proof. Given a feasible solution of (multi-attribute-IRP), we can deduce a feasible solution of (flow-relaxation) by fixing the quantities sent by each depot to each customer per day, delay and commodity. The inventory costs are modelled exactly with (flow-relaxation) thanks to the delayed arcs, thus equal to the ones of (multi-attribute-IRP). The
transportation costs are lower bounded with Equation\:\eqref{eq:fraction_vehicle}. The structure of the graph detailed on Figure\:\ref{fig:commodity_graph}
only allows geographically direct routes between depots and customers, whereas the route constraints allow up to\:$S_{\max}$ stops.
Nonetheless, considering additional arcs\:$c_1 \to c_2$ with fraction costs\:${c^{\texttt{tr}}_{mc_1c_2} = \frac{\ell_m}{L}(c^{\texttt{stop}} + c^{\texttt{km}} \Delta_{c_1c_2})}$
and delay\:$\tau_{c_1c_2}$ in an extended flow graph would also produce solutions with only direct routes. Indeed, by triangular inequality, it would always be cheaper to send commodities through
geographically direct\:$(d \to c)$ (possibly delayed) arcs in this framework of ``vehicle fraction'' costs, rather than sending quantities to intermediate customer\:$c_1$
before reaching the destination\:$c_2$.

\subsection{Bin Packing}\label{flows_to_routes}
We highlight here how the minimum cost flows can be used to derive an IRP solution, a step further from the lower bound computation. The\:$|M|$ minimum cost flow solutions resulting from (flow-relaxation) enable us to set the quantities sent by each depot to each customer per day and commodity, but do not directly lead to a set of routes. Indeed, for now, we do not know how quantities are loaded in various vehicles. We highlight the fact that the delayed arcs are not used to build an initial solution. They are only introduced to compute a lower bound. To deduce a set of direct routes, we approximately solve one bin packing problem per tuple\:${(d,c,t) \in D \times C \times [T]}$, using the first-fit-decreasing heuristic. The instance of the bin packing is given by the set of commodities to be sent on day\:$t$ from depot\:$d$ to customer\:$c$, their respective lengths, and the length of one vehicle. The solution to the bin packing problem leads to a low number of vehicles each of length\:$L$, with corresponding loading made of possibly\:$|M|$ distinct commodities. At this point, we get a set of direct routes as a first feasible solution 
to (multi-attribute-IRP).

\section{Routing Local Search Subroutine}\label{sec:routing_LS}
Our solution processes emphasized on Figure\:\ref{fig:algorithms} rely on  the routing local search subroutine of Algorithm\:\ref{algo:multi_depot_local_search}. We detail here both the local search procedure, and the neighborhoods listed in Table\:\ref{tab:route_moves}.

\subsection{Neighborhoods}
\label{subsec:route_moves}
Before introducing the local search procedure, we focus on the TSP and SDVRP neighborhoods in an IRP framework. The routing local search subroutine combines those neighborhoods.
\paragraph{Routes impact inventories.} We highlight the fact that the\:$12$ neighborhoods detailed below in Table\:\ref{tab:route_moves} alter the routes of a solution and the inventories at the depots or at the customers (contrary to the SDVRP framework). Therefore, whenever a neighborhood is considered, we evaluate the effects on inventories and routes
so as to check feasibility and to estimate the cost change. For instance, the optimal order of a route not only depends on the distances\:$\Delta$, but also on the delays introduced in the inventory dynamics of the customers involved. Therefore, even elementary neighborhoods require calculation. They can be seen as a generalization of 
the TSP and SDVRP concepts to the continuous-time IRP. %Since the computations involved in TSP, single-depot SDVRP 
% and multi-depot SDVRP neighborhoods are substantially distinct, we create different versions of certain neighborhoods.

\begin{table*}[htb]
  \centering
  \caption{Routing local neighborhoods: single-depot and multi-depot variants are considered.}  \label{tab:route_moves}
    \setlength{\tabcolsep} {0.3cm} 
    \scalebox{0.80} 
    { 
    \begin{tabular}{cll} 
  \toprule
  \textbf{Type} & \textbf{Name} & \textbf{Description} \\
  \midrule
   & $\texttt{relocate}$              & change the position of one stop in a route\\
  TSP & $\texttt{swap}$                  & exchange the positions of two stops in a route\\ 
  %  & $\texttt{optimize route}$      & optimize the order of the stops in a route by enumeration.\\
   & $\texttt{2-opt*}$                 & cut a route into three parts and revert the order of the middle one \\
  \midrule
   & $\texttt{insert}$           & give a stop $s$ from route $r_1$ to route $r_2$\\
   & $\texttt{swap single depot}$             &      $\texttt{swap}$ adapted to two routes with same depot\\
  SDVRP single-depot & $\texttt{merge}$                & merge two routes $r_1$ and $r_2$ on the same day\\
  & $\texttt{merge multi day}$       & $\texttt{merge}$ extended to routes with different start dates\\
   & $\texttt{delete route}$         & delete a route\\
   & $\texttt{change day}$                  & move a route in time, one day before or after\\
  \midrule
   & $\texttt{insert multi depot}$            & $\texttt{insert}$ extended to routes with distinct depots \\
  SDVRP multi-depot & $\texttt{swap multi depot}$       & $\texttt{swap}$ extended to routes with distinct depots \\

   & $\texttt{2-opt* multi depot}$                 & cut two routes each into two parts and exchange their end parts \\
  \bottomrule
  \end{tabular} %
  } 
\end{table*}

\subsection{Routing Local Search}
\label{subsec:sdls_mdls}

From the list of neighborhoods emphasized in Table\:\ref{tab:route_moves}, we design the routing local search, Algorithm\:\ref{algo:multi_depot_local_search}. Usually, local search procedures start from the smallest neighborhoods and use larger ones to escape from local minima. Our case is different because of the size and complexity of the problem, and the short time given to solve it. The $\texttt{delete route}$,  $\texttt{change day}$ and TSP neighborhoods (that we name single-route neighborhoods) can be explored exhaustively, until no improvement is found. The SDVRP neighborhoods involve a lot of computations and cannot be fully browsed. Besides, we want each route in the solution after the routing local search procedure to be optimized from a TSP viewpoint.
Therefore, Algorithm\:\ref{algo:multi_depot_local_search} first browses randomly the large SDVRP neighborhoods to create new routes, and then refines thoroughly with the single route neighborhoods. Since the neighborhoods tend to reduce the number of routes, the routing local search algorithm deletes routes at several stages. It only applies feasible moves that improve the cost.

\paragraph{} As said above, some features of this local search differ from the common SDVRP local search algorithms. Because of the computations involved for\:$\texttt{insert}$, \:$\texttt{swap single depot}$, and for the multi-depot neighborhoods, and the number of pairs of routes, the routing local search subroutine only samples a fraction of them. To do so, it explicitly samples a subset of the pairs of routes of the solution according to a uniform distribution on each day. We tune this approach with parameter\:$p$ in Algorithm\:\ref{algo:multi_depot_local_search} to find a good cost gain per CPU time ratio. We do so instead of restricting the routes candidates with geographic criteria, because the inventory costs cannot be neglected. Two routes that visit customers that are far from each other may still be suitable candidates for a\:$\texttt{swap}$ for instance, due to the change in inventory cost. The\:$\texttt{change day}$ function is applied per route one day forward or backward, until no improvement is found.

\begin{algorithm}[!h]
\SetKwFunction{commodityreinsertionMILP}{$\texttt{commodity insertion MILP}$}\SetKwFunction{fillformerroutes}{$\texttt{fill former routes}$}\SetKwFunction{fillnewroutes}{$\texttt{fill new routes}$}\SetKwFunction{singledepotLS}{$\texttt{single depot local search}$}
\SetKwInOut{Input}{input}\SetKwInOut{Output}{output}
\Input{$\mathcal{I}$ an IRP instance,\:$\textbf{r} = (r_k)_{1 \leq k \leq K}$ the current solution with $K \in \mathbb{Z}^+$ routes, $n_{it} \in \mathbb{Z}^+$ a number of iterations, $p$ a percentage.}
\Output{The solution $\textbf{r}$ updated.}
\For{$i = 1:n_{it}$}{
  SDVRP multi-depot neighborhoods over $p\%$ of the pairs of routes at random\;
  $\texttt{delete route}$ per day and depot until no improvement\;
  SDVRP single-depot neighborhoods over $p\%$ of the pairs of routes at random\;
  $\texttt{delete route}$ per day and depot until no improvement\;
  SDVRP single-depot neighborhoods over $p\%$ of the pairs of routes at random\;
  single-route neighborhoods until no improvement;
}
\caption{Routing local search}\label{algo:multi_depot_local_search}
\end{algorithm}

\section{MILP-Based Neighborhoods and Perturbations}\label{sec:MILP_based}
We now introduce three subroutines that are based on optimization problems written as MILPs. They all leverage the commodity graph structure emphasized in Section\:\ref{sec:flow_graphs_formulations}.

\subsection{Reload Fixed-Path Vehicles Subroutine}\label{reload}
\paragraph{Reload neighborhood problem.}
Let us define the problem behind this large neighborhood.
We consider a subset of routes\:$\textbf{r}_{\texttt{reload}}$ of the 
current IRP solution\:$\textbf{r}$, in our case the routes that start from a given depot\:$d$. We solve the following problem: choose the routes to keep in the solution among\:$\textbf{r}_{\texttt{reload}}$, and re-estimate the delivered quantities (for the whole set of commodities) of the routes kept, to minimize the total cost.
In our large neighborhood setting, we fix the remaining routes of the current solution\:$\textbf{r}$.

We denote by\:$x_r$ for\:$r \in \textbf{r}_{\texttt{reload}}$ the indicator variable for 
keeping route\:$r$, and by\:$\bfy = (\bfy_m)_{m \in M}$ the set of commodity flow variables. 
The commodity flow graphs\:$\big(\mathcal{D}^m(\textbf{r}, \textbf{r}_{\texttt{reload}})\big)_{m \in M}$ involved are defined below. We model the problem with the following formulation:
\begin{subequations}
\begin{alignat}{2}
  &\min_{\bfy, \bfx}  & \quad & \phantom{+} \sum_{m \in M} \bfy_m^{\top} \textbf{c}_m + \sum_{r \in \textbf{r}_{\texttt{reload}}}x_r\bigg(c^{\texttt{veh}} + c^{\texttt{stop}}(|P^r|-1) + c^{\texttt{km}} \sum_{a \in A(P^r)}\Delta_a \bigg)  \tag{Reload-MILP}\\*
  & \text{subject to} &       & \bfy_m \in \mathcal{C}\big(\mathcal{D}^m(\textbf{r}, \textbf{r}_{\texttt{reload}}), \bfy_{m}^{\min}(\textbf{r}, \textbf{r}_{\texttt{reload}}), \bfy_{m}^{\max}(\textbf{r}, \textbf{r}_{\texttt{reload}})\big), \quad \forall m \in M\\* %\sum_{a \in \delta^+(v)} \bfy_{ma} = \sum_{a \in \delta^-(v)} \bfy_{ma}, \quad \forall m \in M, \quad \forall v \in \calV^m(\textbf{r}, \textbf{r}_{\texttt{reload}}) \label{eq:cirucaltion_commoditie|P^r|-1eload}\\*
  %&                   &       & \bfy_m^{\min} \leq \bfy_m \leq \bfy_m^{\max}, \quad \forall m \in M \label{eq:capacity_commodity_flow_reload} \\*
  &                   &       &\sum_{m \in M} y_{ma}\ell_m \leq x_rL, \quad \forall a = (d \to (c, r)), \quad \forall r \in \textbf{r}_{\texttt{reload}} \label{eq:capa_new_routes_customer_reinsertion_reload}\\*
  &                   &       &\bfy \in \mathbb{Z} \label{eq:integers_cust_MILP_reload}\\*
  &                   &       &\bfx \in \{0, 1\}   \label{eq:indicator_vehicle_reload}
\end{alignat}
\end{subequations}

This formulation is very close to the generic MILP introduced in Section\:\ref{sec:flow_graphs_formulations}.
The objective function is composed of one flow cost per commodity\:$m \in M$ (inventory and shortage costs),
and of the routing cost of each route kept among\:$\textbf{r}_{\texttt{reload}}$. %The constraints 
%\eqref{eq:cirucaltion_commoditie|P^r|-1eload} and \eqref{eq:capacity_commodity_flow_reload} enforce respectively the circulation 
%and capacity of each commodity flow. 
Constraint \eqref{eq:capa_new_routes_customer_reinsertion_reload} ensures 
that the commodity flows from depots to customers only exist along routes that are kept, 
and that the capacity of the vehicles is respected. The last two constraints define integer and binary variables. We highlight that this 
MILP exactly formulates the reloading of a given subset of routes. Solving the problem (Reload-MILP) leads to a new feasible solution with lower cost. Indeed, we exactly model the IRP constraints and costs, and optimize with respect to a subset of variables with (Reload-MILP): the commodity flows involved in the routes $\textbf{r}_{\texttt{reload}}$, and the use of those routes. The rest of the variables are fixed. The current solution is an admissible solution of this MILP, used for warm-start.

\begin{figure}[htb]
  \centering
  \subfloat[Graph details]{
\begin{minipage}{1.\textwidth}
  \includegraphics[width=\linewidth]{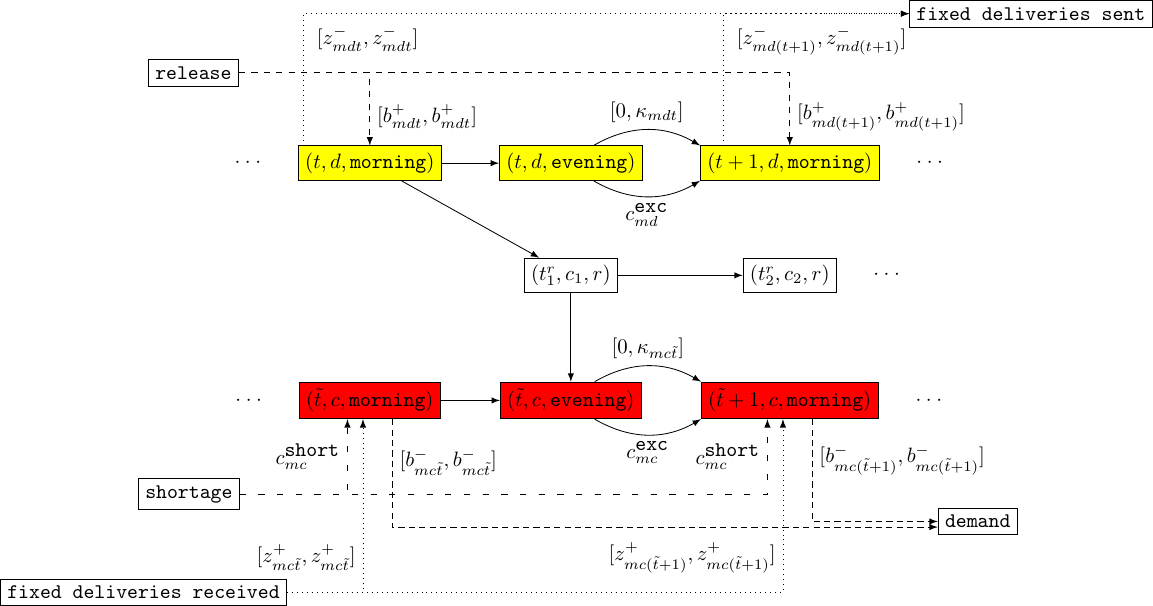}
\end{minipage}
  }
  \newline
  \subfloat[Additional commodity graph arcs compared with Table\:\ref{tab:arcs_generic_subgraphs}]{
\begin{minipage}{1.\textwidth}
  \vspace{0.5cm}
  \setlength{\tabcolsep} {0.3cm} 
  \centering
  \scalebox{0.7} 
  { 
  \begin{tabular}{ccccccc} 
  \toprule
  \textbf{Subgraph}  &  \textbf{Arc description}   & \textbf{Origin}   & \textbf{Destination} & \textbf{Min} & \textbf{Max} & \textbf{Cost} \\
  \midrule
  Routes  & Fixed deliveries sent          & $(t, d, \texttt{morning})$         & $\texttt{fixed deliveries sent}$     & $z^-_{mdt}$  & $z^-_{mdt}$  &                                                       \\ 
  Routes  & Fixed deliveries received      & $\texttt{fixed deliveries received}$ & $(t, c, \texttt{evening})$     & $z^+_{mct}$  & $z^+_{mct}$  &                                                       \\ 
  \midrule 
  Routes  & Transport $d \to (c_1,r)$        & $(t^r, d, \texttt{morning})$   & $(t_1^r, c_1, r)$      &              &              &  \\ 
  Routes  & Transport $(c_s,r) \to (c_{s+1},r)$   & $(t_s^r, c_s, r)$   & $(t_{s+1}^r, c_{s+1}, r)$      &              &              &  \\ 
  Routes  & Transport $(c_s,r) \to c$            & $(t_s^r, c, r)$         & $(t_s^r, c, \texttt{evening})$        &           &              &  \\ 
  \midrule 
  Artificial  & Circulation               & $\texttt{source}$       & $\texttt{fixed deliveries received}$   &              &              &                                  \\ 
  Artificial  & Circulation               & $\texttt{fixed deliveries sent}$       & $\texttt{sink}$   &              &              &                                  \\ 
  \bottomrule
  \end{tabular} %
  }
\end{minipage}
  }
\caption{Details of the commodity graph\:$\mathcal{D}^m(\textbf{r}, \textbf{r}_{\texttt{reload}})$ for the reload fixed-path vehicles neighborhood.}
\label{fig:commodity_graph_refill}
\end{figure}

The commodity graph\:$\mathcal{D}^m(\textbf{r}, \textbf{r}_{\texttt{reload}})$ for\:$m \in M$ depends both on the current solution\:$\textbf{r}$, and on the routes to potentially keep and reload\:$\textbf{r}_{\texttt{reload}}$. 
As previously, the backbone structure is the same as in Section\:\ref{sec:flow_graphs_formulations}: one subgraph per customer, one per depot, one for the routes, and additional vertices to create circulations. 
The special route subgraph is detailed on Figure\:\ref{fig:commodity_graph_refill}. In arc annotations, capacities are given between brackets (e.g. $[0, \kappa_{mdt}]$), and costs without (e.g. $c_{md}^{\texttt{exc}}$). Dotted arrows 
are related to the shared artificial vertices, continuous ones to depots, customers and route subgraphs. When not stated in the table, \textbf{Min} is\:$0$, \textbf{Cost} is\:$0$ and \textbf{Max} is\:$\infty$. We explicitly create individual route paths, with vertices of the form\:$(t_s^r, c_s, r)$. It models the fact that when using route\:$r$, commodities are delivered to customer\:$c_s$ on day\:$t_s^r$ at position\:$s$ of the route. The details of the arcs can be found in the table of Figure\:\ref{fig:commodity_graph_refill}.
Besides, since we fix the quantities sent by the routes in\:$\textbf{r} \backslash \textbf{r}_{\texttt{reload}}$, we need two additional vertices in the route subgraph that we name\:$\texttt{fixed deliveries sent}$, and\:$\texttt{fixed deliveries received}$. The former is connected to depot vertices in order to take other quantities sent into account. The second is connected to customer vertices to model other quantities received. %The main differences are the following:

%\begin{itemize}
%  \item the customers and depots vertices are restricted to those present in the routes $\textbf{r}_{\texttt{reload}}$.
%  \item two additional artificial vertices are added: respectively for the fixed deliveries to be sent and received, 
%  $\texttt{fixed deliveries sent}$ and $\texttt{fixed deliveries received}$. The latter 
%  are due to the existence of deliveries through routes that belong to $\textbf{r} \backslash \textbf{r}_{\texttt{reload}}$.
%  \item for each route $r \in \textbf{r}_{\texttt{reload}}$, one vertex $(t_c, c, r)$ is added per stop (involving customer $c$ visited on $t_c$ by $r$) 
%  of the route to represent a path.
%  \item no direct arc $(t,d,\texttt{morning}) \to (t',c,\texttt{evening})$ is considered.
%  \item the route path is created linking route vertices in the proper order. Given $r \in \textbf{r}_{\texttt{reload}}$,
%  we start from $(t_r, d_r, \texttt{morning})$ and visit the vertices $(t_c, c, r)$ in the stops order. 
%  \item we create an arc from each route vertex $(t_c,c,r)$ to the evening vertex of the customer $(t_c,c,\texttt{evening})$.
%  \item we add one arc from each depot morning vertex to $\texttt{fixed deliveries sent}$ and from $\texttt{fixed deliveries received}$
%  to each customer evening.
%\end{itemize}

\subsection{Customer Reinsertion Subroutine}\label{customer_reinsertion}

As mentioned in the overview Section\:\ref{sec:algos_overview},
our neighborhoods are based on a decomposition of the IRP along its 
main axes. Previous sections focus on the routes. Here, we design a perturbation based on the customers. We call it perturbation because it may slightly increase the cost of the IRP solution. 

\paragraph{Customer reinsertion problem.} Let us define the customer reinsertion problem. Once the customer\:$c \in C$ is removed from the solution -- that is to say, removed from the routes that deliver to it, leading to zero delivery\:$\bfz^+$ in the inventory dynamics \eqref{eq:durations_inventory_dynamics_customer} -- we need to reinsert it in the solution, using only former routes and new direct routes. This means choosing: 1)\:The insertion position of customer\:$c$ in each route of the solution in which it is inserted, keeping the relative order of the other stops unchanged. 2)\:The quantity of each commodity to be delivered by those former routes where\:$c$ is inserted. 3)\:New direct routes (path, timing and quantities) to deliver\:$c$. It can be formulated as a MILP akin to the generic one defined in Section\:\ref{sec:flow_graphs_formulations}. The reload fixed-path vehicles and customer reinsertion formulations are different. In the former, we consider only a set of given fixed-path routes in the route subgraph, the associated customers, and one depot. The decisions are selecting or not each route and fixing the quantities to be delivered to each stop. Both inventory and routing costs induced are exactly modelled. In the latter, a customer is removed from every delivery of the current solution. With the reinsertion MILP, we decide in which former route we insert the customer, at which position, and how much of each commodity we deliver to it. We also decide if we create new direct routes from depots to this particular customer and how much of each commodity we send through them. New direct routes are approximated by one large bin per depot origin and departure date. Inventory dynamics are exactly modelled, at the customers impacted by the insertion due to delays, and at the depots. 

% \begin{algorithm}[htb]
% \SetKwData{Action}{action}
% \SetKwFunction{customerinsertionmilp}{$\texttt{customer insertion MILP}$}\SetKwFunction{fillformerroutes}{$\texttt{fill former routes}$}\SetKwFunction{fillnewroutes}{$\texttt{fill new routes}$}\SetKwFunction{removecustomer}{$\texttt{remove customer}$}
% \SetKwInOut{Input}{input}\SetKwInOut{Output}{output}
% \Input{$\mathcal{I}$ an IRP instance, $\textbf{r} = (r_k)_{1 \leq k \leq K}$ the current solution with $K \in \mathbb{Z}^+$ routes, $c \in C$ a customer.}
% \Output{The solution $\textbf{r}$ updated.}
% $\textbf{r}$ = \removecustomer($\mathcal{I}$, $\textbf{r}$, $c$)\;
% $\bfy$ = \customerinsertionmilp($\mathcal{I}$, $\textbf{r}$, $c$)\;
% $\textbf{r}$ = \fillformerroutes($\mathcal{I}$, $\textbf{r}$, $c$, $\bfy$)\;
% $\textbf{r}$ = \fillnewroutes($\mathcal{I}$, $\textbf{r}$, $c$, $\bfy$)\;
% \caption{One step of customer reinsertion}\label{algo:customer_reinsertion}
% \end{algorithm}

\paragraph{The customer insertion MILP.} Let\:$\bfy = (\bfy_m)_{m \in M_c}$ be the commodity flow variable. Instead of the indicator route variables in Section\:\ref{sec:flow_graphs_formulations}, we use another type of graph to model the vehicles with a flow variable\:$\bfx$. Indeed, a flow is a convenient tool to model the fact that we have multiple insertion positions in a given route for customer\:$c$, and we can choose at most one of them. We link the flow variables on the route subgraphs of the commodity graphs with this vehicle flow. It leads to the MILP:
\begin{subequations}
  \begin{alignat}{2}
  &\min_{\bfy, \bfx}  & \quad & \phantom{+} \sum_{m \in M_c} \bfy_m^{\top} \textbf{c}_m + c_x^{\top} \bfx  \tag{Cust-MILP}\\*
  & \text{subject to} &       & \bfy_m \in \mathcal{C}(\mathcal{D}^m_c, \bfy_{m}^{\min}, \bfy_{m}^{\max}), \quad \forall m \in M_c\\* 
  &                   &       & \bfx \in \mathcal{C}(\mathcal{D}^{\texttt{veh}}_c, \bfx^{\min}, \bfx^{\max})\\*%\sum_{a \in \delta^+(v)} \bfy_{ma} = \sum_{a \in \delta^-(v)} \bfy_{ma}, \quad \forall m \in M, \quad \forall v \in \calV^m(\textbf{r}, \textbf{r}_{\texttt{reload}}) \label{eq:cirucaltion_commoditie|P^r|-1eload}\\*
  %& \text{subject to} &       & \sum_{a \in \delta^+(v)} \bfy_{ma} = \sum_{a \in \delta^-(v)} \bfy_{ma}, \quad \forall m \in M_c, \quad \forall v \in \calV^m_c \label{eq:cirucaltion_commodities}\\*
  %&                   &       &\sum_{a \in \delta^+(v)} \bfx_{a} = \sum_{a \in \delta^-(v)} \bfx_{a}, \quad \forall v \in \calV^{\texttt{veh}}_c  \label{eq:circulation_vehicles} \\*
  %&                   &       &\bfy_m^{\min} \leq \bfy_m \leq \bfy_m^{\max}, \quad \forall m \in M_c \label{eq:capacity_commodity_flow} \\*
  %&                   &       &\bfx^{\min} \leq \bfx \leq \bfx^{\max} \label{eq:capacity_vehicle_flow}\\* 
  &                   &       &\sum_{m \in M_c} y_{ma}\ell_m \leq x_aL, \quad \forall a = (d \to c), \quad \forall d \in D \label{eq:capa_new_routes_customer_reinsertion}\\*
  &                   &       &\sum_{m \in M_c} y_{ma}\ell_m \leq x_aL^r_{\texttt{free}}, \quad \forall a = (r \to (r,s)), \quad \forall r \in \textbf{r}, \quad \forall s \in [|P^r|] \label{eq:capa_former_routes_customer_reinsertion}\\*
  &                   &       &\bfx \in \mathbb{Z}, ~ \bfy \in \mathbb{Z} \label{eq:integers_cust_MILP}
\end{alignat}
\end{subequations}

\paragraph{}

%Equations \eqref{eq:cirucaltion_commodities} and \eqref{eq:circulation_vehicles} are the circulation constraints, 
%and \eqref{eq:capacity_commodity_flow} and \eqref{eq:capacity_vehicle_flow} the capacity constraints.
The vehicle flow as well as each commodity flow must respect circulation constraints. In Equation \eqref{eq:capa_new_routes_customer_reinsertion}, we force the amount of commodities to be sent through new direct routes from the depots to the customer\:$c$ not to exceed the total content size of the vehicles involved. Indeed, we consider the total vehicle capacity with this constraint, and not individual vehicles each of capacity $L$. Similarly, Equation \eqref{eq:capa_former_routes_customer_reinsertion} does so for the former routes\:$r \in \textbf{r}$ at position\:$s$ having remaining loading space\:$L^r_{\texttt{free}}$.

\begin{figure}[htb]
  \centering
  \subfloat[Graph details]{
\begin{minipage}{1.\textwidth}
  \includegraphics[width=\linewidth]{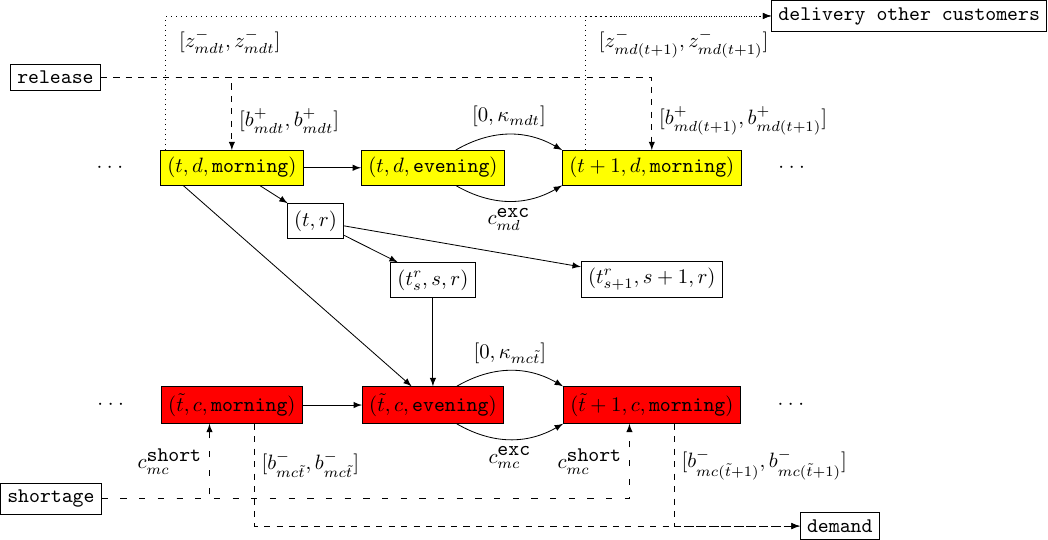}
\end{minipage}
  }
  \newline
  \subfloat[Additional commodity graph arcs compared with Table\:\ref{tab:arcs_generic_subgraphs}]{
\begin{minipage}{1.\textwidth}
  \centering
  \vspace{0.5cm}
  \setlength{\tabcolsep} {0.3cm} 
  \scalebox{0.70} 
  { 
  \begin{tabular}{ccccccc} 
  \toprule
  \textbf{Subgraph} & \textbf{Arc description}     & \textbf{Origin}   & \textbf{Destination} & \textbf{Min} & \textbf{Max} & \textbf{Cost} \\
  \midrule
  Routes & Delivery other customers       & $(t, d, \texttt{morning})$         & $\texttt{delivery other customers}$     & $z^-_{mdt}$  & $z^-_{mdt}$  &                                                       \\ 
  \midrule 
  Routes & Transport $d \to c$ new route            & $(t, d, \texttt{morning})$         & $(t+\floor{\frac{\tau_{dc}}{\tau_{\max}}}, c, \texttt{evening})$      &              &              &  \\ 
  Routes & Transport $d \to r$ former route            & $(t, d, \texttt{morning})$         & $(t, r)$      &              &              &  \\ 
  Routes & Transport $r \to (r,s)$ former route       & $(t, r)$                           & $(t_s^r, s, r)$      &              &              &  \\ 
  Routes & Transport $(r,s) \to c$ former route         & $(t_s^r, s, r)$         & $(t_s^r, c, \texttt{evening})$        &           &              &  \\ 
  \midrule 
  Artificial & Circulation               & $\texttt{delivery other customers}$       & $\texttt{sink}$   &              &              &                                  \\ 
  \bottomrule
  \end{tabular} %
  }
\end{minipage}
  }
\caption{Details of the commodity graph\:$\mathcal{D}^m_c$ for the customer reinsertion neighborhood.}
\label{fig:commodity_graph_customer_reinsertion}
\end{figure}

\paragraph{Details of the commodity flow graphs.} Based on the generic graph structure Section\:\ref{sec:flow_graphs_formulations}, we define new commodity graphs\:$\mathcal{D}^m_c = (\calV^m_c, \calA^m_c)$ for each commodity\:$m \in M_c$ for the reinsertion of customer\:$c \in C$. We highlight the fact that here only one customer is involved in the problem, so only one customer subgraph is present. The depots and customer subgraphs are introduced in Section\:\ref{sec:flow_graphs_formulations}. We now define the route subgraph specific to this customer reinsertion MILP. First, we have a vertex\:$\texttt{delivery other customers}$ in the route subgraph. It is connected to the depot subgraphs to model the quantities sent to other customers. Figure\:\ref{fig:commodity_graph_customer_reinsertion} shows the specific commodity graph. In arc annotations, capacities are given between brackets (e.g. $[0, \kappa_{mdt}]$), and costs without (e.g. $c_{md}^{\texttt{exc}}$). Dotted arrows are related to the shared artificial vertices, continuous ones to depots, customers and route subgraphs. When not stated in the table, \textbf{Min} is\:$0$, \textbf{Cost} is\:$0$ and \textbf{Max} is\:$\infty$. In this commodity graph, we consider two types of routes. 1)\:new direct routes are modelled by direct arcs of the type\:${(t, d, \texttt{morning}) \to (t', c, \texttt{evening})}$ and do not involve additional vertices. 2)\:former routes in which we can insert customer\:$c$. They are modelled with one vertex\:$(t, r)$ for route\:$r$ starting on day\:$t$, connected to the starting depot morning vertex and to each vertex of the form\:$(t_s^r,s,r)$. The vertex\:$(t_s^r,s,r)$ is related to the possible insertion position\:$s$ in route\:$r$ leading to an arrival day\:$t_s^r$ at customer\:$c$. It is connected to the customer subgraph. The details of the route arcs are in the table of Figure\:\ref{fig:commodity_graph_customer_reinsertion}. We cannot only add one vertex per former route, since the optimal insertion position depends on the commodity flows, and not only on routing costs. Therefore, the optimal insertion position cannot be pre-computed. Besides, we do not introduce an approximate cost for the transportation arcs. Instead, we define another graph\:$\mathcal{D}^{\texttt{veh}}_c = (\calV^{\texttt{veh}}_c, \calA^{\texttt{veh}}_c)$ as follows. 

\paragraph{Details of the vehicle flow graph.}
We consider the\:$\calV^{\texttt{veh}}_c$ vertices:

\begin{itemize}
  \item Artificial vertices labelled:\:$\texttt{source}$ and\:$\texttt{sink}$.
  \item Depots vertices labelled: $(t,d, \texttt{morning})$ for\:$t \in [T]$, $d \in D$.
  \item Customer vertices labelled: $(t,c,\texttt{evening})$ for\:$t \in [T]$.
  \item Routes vertices labelled: $(t,r)$ for each former route\:$r \in \textbf{r}$ that 
  can include an additional stop and that does not reach the vehicle capacity\:$L$ already. 
  \item Nodes per insertion position in former routes:\:$(t_s^r,s,r)$ for each position\:$s$ at which we can insert 
  customer\:$c$ in\:$r$ without exceeding the time horizon\:$T$. The date\:$t_s^r$ on which the route delivers customer\:$c$
  can be pre-computed. 
\end{itemize}

\begin{table*}[htb]
  \centering
  \caption{Arcs of the customer reinsertion vehicle flow graph\:$\mathcal{D}^{\texttt{veh}}_c$. \\ When not stated, \textbf{Min} is\:$0$, \textbf{Cost} is\:$0$ and 
  \textbf{Max} is\:$\infty$.}    \label{tab:arcs_vehicle_graph_customer_reinsertion}
  \setlength{\tabcolsep} {0.3cm} 
  \scalebox{0.8} 
  { 
  \begin{tabular} {cccccc} 
  \toprule
  \textbf{Arc description}     & \textbf{Origin}   & \textbf{Destination} & \textbf{Min} & \textbf{Max} & \textbf{Cost}  \\
  \midrule 
  Start routes                 & $\texttt{source}$                     & $(t, d, \texttt{morning})$     &   &   &                                                       \\ 
  End routes                   & $(t, c, \texttt{evening})$           & $\texttt{sink}$       &   &   &                                                       \\ 
  \midrule
  Transport $(d \to c)$ new route        & $(t, d, \texttt{morning})$           & $(t+\floor{\frac{\tau_{dc}}{\tau_{\max}}}, c, \texttt{evening})$     &   &   &   $c^{\texttt{veh}} + c^{\texttt{stop}} + c^{\texttt{km}} \Delta_{dc} $      \\ 
  Transport $(d \to r)$ former route       & $(t, d, \texttt{morning})$           & $(t, r)$     &   & 1  &        \\ 
  Transport $(r \to (r,s))$ former route    & $(t, r)$           & $(t_s^r, s, r)$     &   &   &  $c_{sc}^r$    \\ 
  Transport $((r,s) \to c)$ former route       & $(t_s^r, s, r)$                             & $(t_s^r, c, \texttt{evening})$     &   &   &      \\ 
  \midrule
  Circulation                   & $\texttt{source}$           & $\texttt{sink}$       &   &   &                                                       \\ 
  Circulation                   & $\texttt{sink}$           & $\texttt{source}$       &   &   &                                                       \\ 
  \bottomrule
  \end{tabular} %
  } 
  \end{table*}

The arcs of the customer reinsertion vehicle flow graph are defined in Table\:\ref{tab:arcs_vehicle_graph_customer_reinsertion}. We see that most of its structure is shared with the commodity flow graph described above. Besides, the costs~$c_{sc}^r$ are the ones induced by the insertion of customer\:$c$ at position\:$s$ in the stops of route~$r$. They involve routing and inventory considerations at the other customers delivered by route\:$r$, because of the delays. They can be exactly computed considering the former list of\:$r$'s stops, and enumerating the insertion possibilities. We emphasize the\:$1$-maximum capacity on the arcs of the form\:${(t, d, \texttt{morning}) \to (t,r)}$ combined with the circulation constraint enforce that at most one insertion position is chosen in former routes. 

\begin{proposition}
  The problem (Cust-MILP) is a relaxation of the optimal customer reinsertion problem in the IRP solution,
  where the new direct routes are aggregated per routing arc\:$(d \to c)$.
\end{proposition}

We sketch the proof. The constraints for new direct routes do not exactly model the routing structure: instead of a set of vehicles having each a content size\:$L$, it is as if we had one large vehicle with the total content size on each arc\:$d \to c$. More precisely, given\:$d \in D$, filling\:$p \in \mathbb{Z}^+$ vehicles of content size\:$L$ and cost\:$c^{\texttt{veh}} + c^{\texttt{stop}} + c^{\texttt{km}} \Delta_{dc}$ each enables less loading freedom than filling one large vehicle of content size\:$pL$ with same total cost\:$p(c^{\texttt{veh}} + c^{\texttt{stop}} + c^{\texttt{km}} \Delta_{dc})$. In the objective function, the costs\:$\textbf{c}_m$ for\:$m \in M_c$ and\:$c_x$ stem from the arc features stated in Table\:\ref{tab:arcs_generic_subgraphs}, the table of Figure\:\ref{fig:commodity_graph_customer_reinsertion}
and in Table\:\ref{tab:arcs_vehicle_graph_customer_reinsertion}. Because of their values and because we relax the routing structure, this problem is a relaxation of the customer reinsertion in the current solution of the IRP.

\paragraph{Rebuilding routes.} All the decisions we make are encapsulated in the commodity flow variable\:$\bfy$. From\:$\bfy$, we can easily update the solution\:$\textbf{r}$, filling former routes with the indicated quantities. Then, a bin packing problem is solved approximately per depot and per day to create the new direct routes from the flows on direct routes arcs. Aggregating per\:$(d \to c)$ arc enables us to derive a MILP of reasonable size, but this relaxation can lead to a potential cost increase. Our LNS accepts the new IRP solution if its true cost is not more than $1\%$ greater than the cost of the previous solution.

\subsection{Commodity Reinsertion Subroutine}\label{commodity_reinsertion}

The last subroutine we detail here is the commodity reinsertion. As for the customer reinsertion, the idea is to perturb the solution at a broad scale, possibly increasing the cost.

% \begin{algorithm}[htb]
% \SetKwData{Action}{action}
% \SetKwFunction{commodityreinsertionMILP}{$\texttt{commodity insertion MILP}$}\SetKwFunction{fillformerroutes}{$\texttt{fill former routes}$}\SetKwFunction{fillnewroutes}{$\texttt{fill new routes}$}\SetKwFunction{singledepotLS}{$\texttt{single depot local search}$}\SetKwFunction{removecommodity}{$\texttt{remove commodity}$}
% \SetKwInOut{Input}{input}\SetKwInOut{Output}{output}
% \Input{$\mathcal{I}$ an IRP instance,\:$\textbf{r} = (r_k)_{1 \leq k \leq K}$ the current solution with\:$K \in \mathbb{Z}^+$ routes,\:$m \in M$ a commodity.}
% \Output{The solution\:$\textbf{r}$ updated.}
%\:$\textbf{r}$ = \removecommodity($\mathcal{I}$,\:$\textbf{r}$,\:$m$)\;
%\:$\bfy$ = \commodityreinsertionMILP($\mathcal{I}$,\:$\textbf{r}$,\:$m$)\;
%\:$\textbf{r}$ = \fillformerroutes($\mathcal{I}$,\:$\textbf{r}$,\:$m$,\:$\bfy$)\;
%\:$\textbf{r}$ = \fillnewroutes($\mathcal{I}$,\:$\textbf{r}$,\:$m$,\:$\bfy$)\;
%\:$\textbf{r}$ = \singledepotLS($\mathcal{I}$,\:$\textbf{r}$)\;
% \caption{One step of commodity reinsertion}\label{algo:commodity_reinsertion}
% \end{algorithm}

\paragraph{Commodity reinsertion problem.}
We are interested in the following problem. Once a commodity\:$m \in M$ is removed from every delivery of the current solution, we want to choose the quantity of\:$m$ to send through each former route, and through new direct routes. Former routes are the routes of the current solution with the commodity\:$m$ removed that have a remaining load. The restriction to new direct routes is a choice to quickly compute a solution. 

\paragraph{The commodity insertion MILP.} As for the customer reinsertion, we couple a commodity flow variable\:$\bfy$ with a vehicle flow\:$\bfx$. The resulting MILP is given by:
\begin{subequations}
\begin{alignat}{2}  
& \min_{\bfy, \bfx} {} & \quad & \phantom{+} \bfy_m^{\top} \textbf{c}_m + c_x^{\top} \bfx \tag{Comm-MILP}\\*
& \text{subject to}    &       & \bfy_m \in \mathcal{C}(\mathcal{D}^m_m, \bfy_{m}^{\min}, \bfy_{m}^{\max})\\*
&                      &       & \bfx \in \mathcal{C}(\mathcal{D}^{\texttt{veh}}_m, \bfx^{\min}, \bfx^{\max})\\*
&                      &       & y_{ma} \leq \left\lfloor{\frac{L}{\ell_m}}\right\rfloor x_a, \quad \forall a = (d \to c) \quad \forall d \in D_m \quad \forall c \in C_m\label{eq:commodity_rein_capa_large_container}\\*
&                      &       & \bfx \in \mathbb{Z}, ~ \bfy \in \mathbb{Z}
\end{alignat}
\end{subequations}
We highlight only one commodity flow and one vehicle flow are involved here. Based on the subgraphs introduced in Section\:\ref{sec:flow_graphs_formulations}, the inventory constraints and costs are exactly modelled. The routing structure is approximated, with a cost paid per unit of a single ``large vehicle with total content size'' per arc\:$(d \to c)$ for\:$c \in C_m$ and\:$d \in D_m$, modelled by constraint\:\eqref{eq:commodity_rein_capa_large_container}.

\paragraph{Details of the commodity flow graph.} Let\:$m \in M$ be the commodity we consider. As previously, we build a commodity graph\:$\mathcal{D}^m_m = (\calV^m_m, \calA^m_m)$. It shares the backbone structure with customer and depot subgraphs as well as artificial vertices defined in Section\:\ref{sec:flow_graphs_formulations}. We now detail the specific route subgraph for the commodity reinsertion problem. It aims at modelling both former routes and new direct routes. For a former route\:$r \in \textbf{r}$, for\:$c$ visited by\:$r$, we denote by \:$s \in [|P^r|-1]$ its position in the list of visited stops, and by\:$t^r_s$ the date on which it is delivered. We have a vertex\:$(t_s^r,c, r)$ in the route subgraph. We then have the arcs of Table\:\ref{tab:additional_arcs_commodity_graph_commodity_reinsertion} to explicitly model the flow of the commodity\:$m$ through former and new routes, in addition to those defined in Table\:\ref{tab:arcs_generic_subgraphs}.

\begin{table*}[htb]
  \centering
  \caption{Additional arcs of the commodity reinsertion commodity flow graph compared to Table\:\ref{tab:arcs_generic_subgraphs}. \\ When not stated, \textbf{Min} is\:$0$, \textbf{Cost} is\:$0$ and 
  \textbf{Max} is\:$\infty$.}    \label{tab:additional_arcs_commodity_graph_commodity_reinsertion}
  \setlength{\tabcolsep} {0.3cm} 
  \scalebox{0.75} 
  { 
  \begin{tabular} {ccccccc} 
  \toprule
  \textbf{Subgraph} & \textbf{Arc description}     & \textbf{Origin}   & \textbf{Destination} & \textbf{Min} & \textbf{Max} & \textbf{Cost}  \\
  \midrule 
  Routes & Transport $d \to c$ new route  & $(t, d, \texttt{morning})$         & $(t+\floor{\frac{\tau_{dc}}{\tau_{\max}}}, c, \texttt{evening})$      &              &              &  \\
  Routes & Transport $(d \to c_{1})$ former route $r$       & $(t^r, d, \texttt{morning})$       & $(t_1^{r}, c_{1}, r)$     &   & $\left\lfloor\frac{L^{r}_{\texttt{free}}}{\ell_m}\right\rfloor$ &                                                       \\ 
  Routes & Transport $(c_{s} \to c_{s+1})$ former route $r$        & $(t_s^{r}, c_{s}, r)$      & $(t_{s+1}^{r}, c_{s+1}, r)$     &   &   &                                                       \\ 
  Routes & Deliver $c_{s}$ former route $r$       & $(t_s^{r}, c_{s}, r)$      & $(t_s^{r}, c_{s}, \texttt{evening})$     &   &   &                                                       \\ 
  \bottomrule
  \end{tabular} %
  } 
  \end{table*}

\paragraph{Details of the vehicle flow graph.} This flow graph is only used to model new direct routes. Indeed, all former routes are reused in this problem: we only choose the commodity flow through them. The only routing decision is related to the creation of new direct routes. Therefore, the \emph{commodity reinsertion vehicles graph vertices}\:$\calV^{\texttt{veh}}_m$ are the following: 1)\:Artificial vertices labelled:\:$\texttt{source}$ and\:$\texttt{sink}$. 2)Depot vertices labelled:\:$(t,d, \texttt{morning})$ for\:$t \in [T]$, $d \in D_m$. 3)\:Customer vertices labelled:\:$(t,c,\texttt{evening})$ for\:$t \in [T]$, $c \in C_m$. Let\:$\calA^{\texttt{veh}}_m$ be the arcs of the vehicle flow graph. They are detailed in Table\:\ref{tab:arcs_vehicle_graph_commodity_reinsertion}.

\begin{table*}[htb]
  \centering
  \caption{Arcs of the commodity reinsertion vehicle flow graph\:$\mathcal{D}^{\texttt{veh}}_m$. \\ When not stated, \textbf{Min} is\:$0$, \textbf{Cost} is\:$0$ and 
  \textbf{Max} is\:$\infty$.}      \label{tab:arcs_vehicle_graph_commodity_reinsertion}
  \setlength{\tabcolsep} {0.3cm} 
  \scalebox{0.92} 
  { 
  \begin{tabular} {cccccc} 
  \toprule
  \textbf{Arc description}     & \textbf{Origin}   & \textbf{Destination} & \textbf{Min} & \textbf{Max} & \textbf{Cost}  \\
  \midrule 
  Start routes                 & $\texttt{source}$                     & $(t, d, \texttt{morning})$     &   &   &                                                       \\ 
  Transport $(d \to c)$        & $(t, d, \texttt{morning})$           & $(t+\floor{\frac{\tau_{dc}}{\tau_{\max}}}, c, \texttt{evening})$     &   &   &   $c^{\texttt{veh}} + c^{\texttt{stop}} + c^{\texttt{km}} \Delta_{dc} $      \\ 
  End routes                   & $(t, c, \texttt{evening})$           & $\texttt{sink}$       &   &   &                                                       \\ 
  \midrule
  Circulation                   & $\texttt{source}$           & $\texttt{sink}$       &   &   &                                                       \\ 
  Circulation                   & $\texttt{sink}$           & $\texttt{source}$       &   &   &                                                       \\ 
  \bottomrule
  \end{tabular} %
  } 
  \end{table*}

\paragraph{Rebuilding routes.} The commodity reinsertion subroutine proceeds as follows. From the information written in the\:$\bfy$ flow variable, it fills the former routes by decoding the corresponding flow. Then, it iteratively creates new direct routes to zero the corresponding commodity flow in\:$\bfy$. In contrast with the previous section, restricting ourselves to new direct routes is a limitation. It applies a\:$\texttt{routing local search}$ to address the potential cost rise induced. The new solution is accepted by the LNS, even if its cost is higher.

\section{Computational Experiments}\label{sec:numerical_results}
First, we want to analyze the quality of the solutions of our LNS. We compare them with the route-based matheuristic in terms of costs. We then perform additional experiments to understand the contribution of the distinct components of our LNS to its performance.

As previously mentioned, all of our code is available on GitHub \parencite{bouvier_louis_2023_8179161}.
The library of instances \parencite{louis_bouvier_2023_8177237} and the solutions we obtained \parencite{louis_bouvier_2023_8177271} can be found on Zenodo.
Note that both the instances and the solutions are downloaded automatically by our test suite, which is run online following every commit.
This is meant to improve the reproducibility and reliability of our experiments.

\subsection{Experimental Setting}\label{subsec:experimental_setting}

\begin{figure}[htb]
  \begin{minipage}{0.33\textwidth}
    \centering
    \includegraphics[width=.85\linewidth]{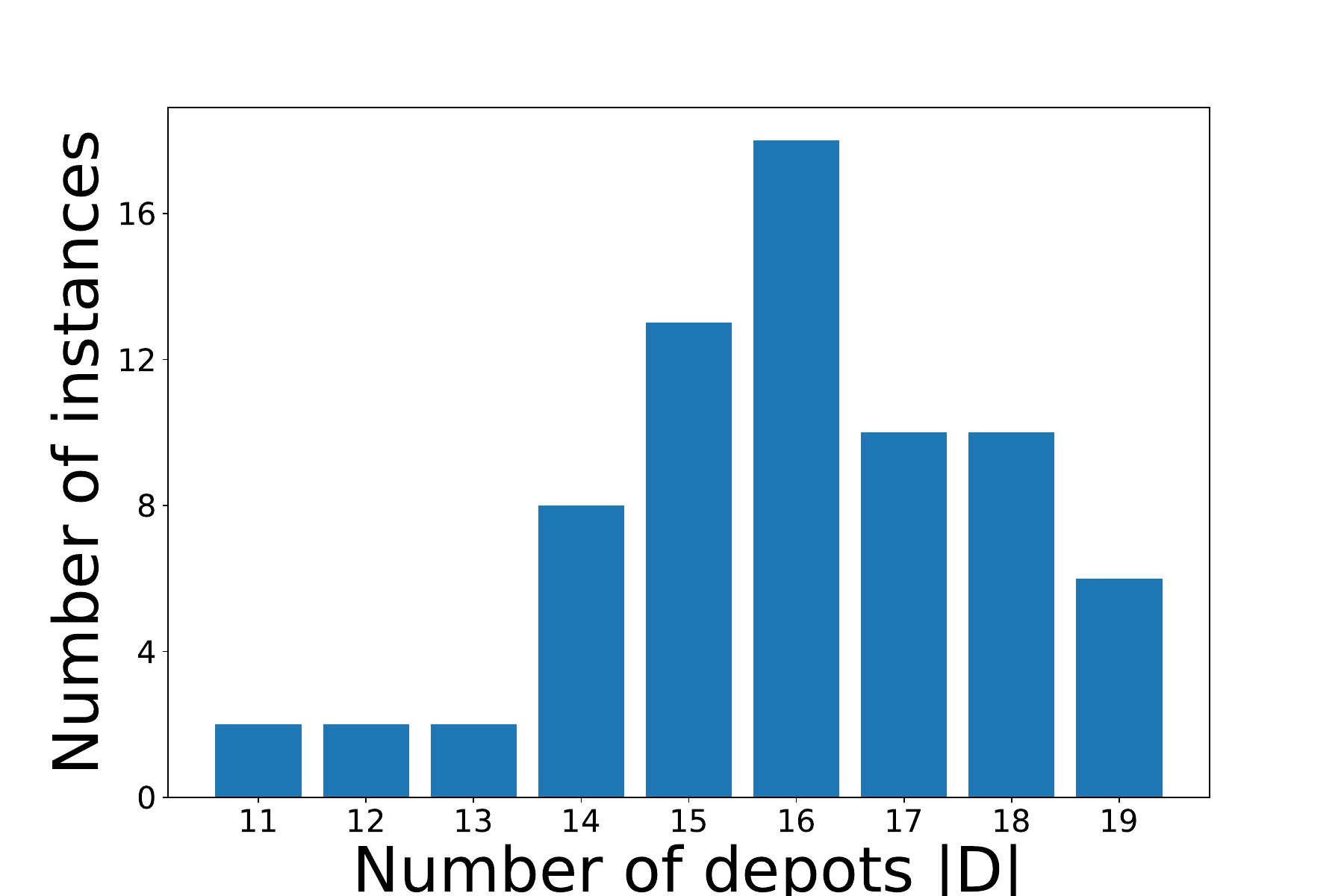}
  \end{minipage}
     \begin{minipage}{0.33\textwidth}
    \centering
    \includegraphics[width=.85\linewidth]{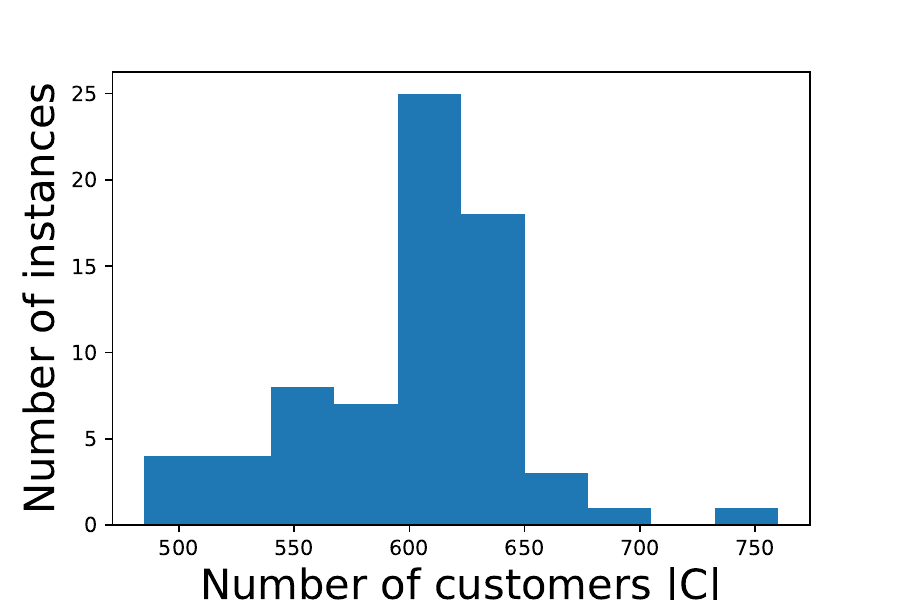}
  \end{minipage}
  \begin{minipage}{0.33\textwidth}
    \centering
    \includegraphics[width=.85\linewidth]{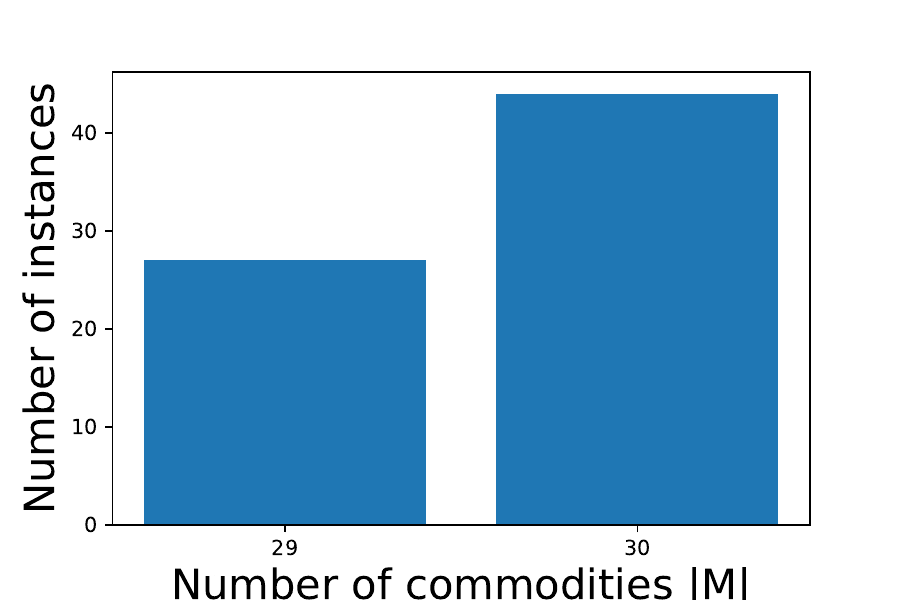}
  \end{minipage}
  \caption{Histograms of the dimensions of the extracted instances.}
  \label{fig:dimensions_instances}
\end{figure}

\paragraph{Instances.} Renault gives us access to\:$71$ IRP instances, each of them corresponding to their return logistics at the European scale and over roughly\:$20$ days. They correspond to the same industrial use case but to different periods.
We show the instances' number of depots, number of customers and number of commodities on Figure~\ref{fig:dimensions_instances}. The maximum number of stops is \:$S_{max} = 3$ to comply with the car manufacturer's requirements. Additional details on the instances are available in Figure\:\ref{fig:instance} in Appendix\:\ref{app:instance}. They highlight why our instances are hard to solve, showing how release and demand are spread respectively among depots and customers, as well as how grouping customers in routes with multiple stops is essential. 

\paragraph{Implementation.} Our code is in the Julia language \parencite{Julia-2017}. We use Gurobi optimizer \parencite{gurobi} to solve LPs and MILPs, JuMP \parencite{DunningHuchetteLubin2017} to model mathematical programs, and $\texttt{Graphs.jl}$ \parencite{Graphs2021} to define graph structures. We proceed to a warm 
start using the current solution for each MILP.
Our experiments are run on a computing cluster with $189$Go of RAM and $32$ processors. Multithreading is enabled for Gurobi only.

\paragraph{Hyperparameters.} Several hyperparameters are introduced for our LNS: the number of routing local search iterations per LNS step, the percentage of random pairs of routes considered in the routing local search, the number of reload fixed-path vehicles, customer, and commodity reinsertion iterations per outer LNS step. To highlight the effect of our new perturbations, we tune the number of customer and commodity reinsertion steps per outer LNS step, as shown in Appendix\:\ref{app:hyperparameters} in Table\:\ref{tab:gridsearch}. The best configuration we find for the short $90$ minutes runs is $|C|$ customer reinsertion (the whole set of customers) and $10$ commodity reinsertion steps. We browse the whole set of depots in the reload fixed-path vehicles subroutine, and proceed to one iteration of routing local search per outer LNS step. Besides, we run our code for $90$ minutes on each instance. We include longer run results with $300$ minutes to solve instances with $S_{\max} = 10$ in the extended experiments of Section\:\ref{subsubsec:res_longer_routes}.

\paragraph{Experimental design.} We run our LNS with the hyperparameters described above on the $71$ real instances. We compare solutions with the route-based matheuristic, as well as with some versions of the LNS where a subroutine is removed in an ablation study. Every algorithm has the same time budget of $90$ minutes. Nonetheless, since both the initialization + local search and route-based matheuristic are constructive, they end before the $90$ minutes time budget. In Appendix\:\ref{app:hyperparameters}, when tuning the LNS hyperparameters, we show in Table\:\ref{tab:gridsearch} a version of the LNS with zero customer and commodity reinsertion step. This version of the LNS can be seen as an iterative extension of the route-based matheuristic. Our aim is to analyze the performance of the LNS, as well as the effect of each of its components. Longer runs with $S_{\max} = 10$ aim at showing rooms for improvement.

\subsection{Performance analysis}
\label{subsec:cost_results}

\begin{figure}[htb]
    \includegraphics[width=.9\linewidth]{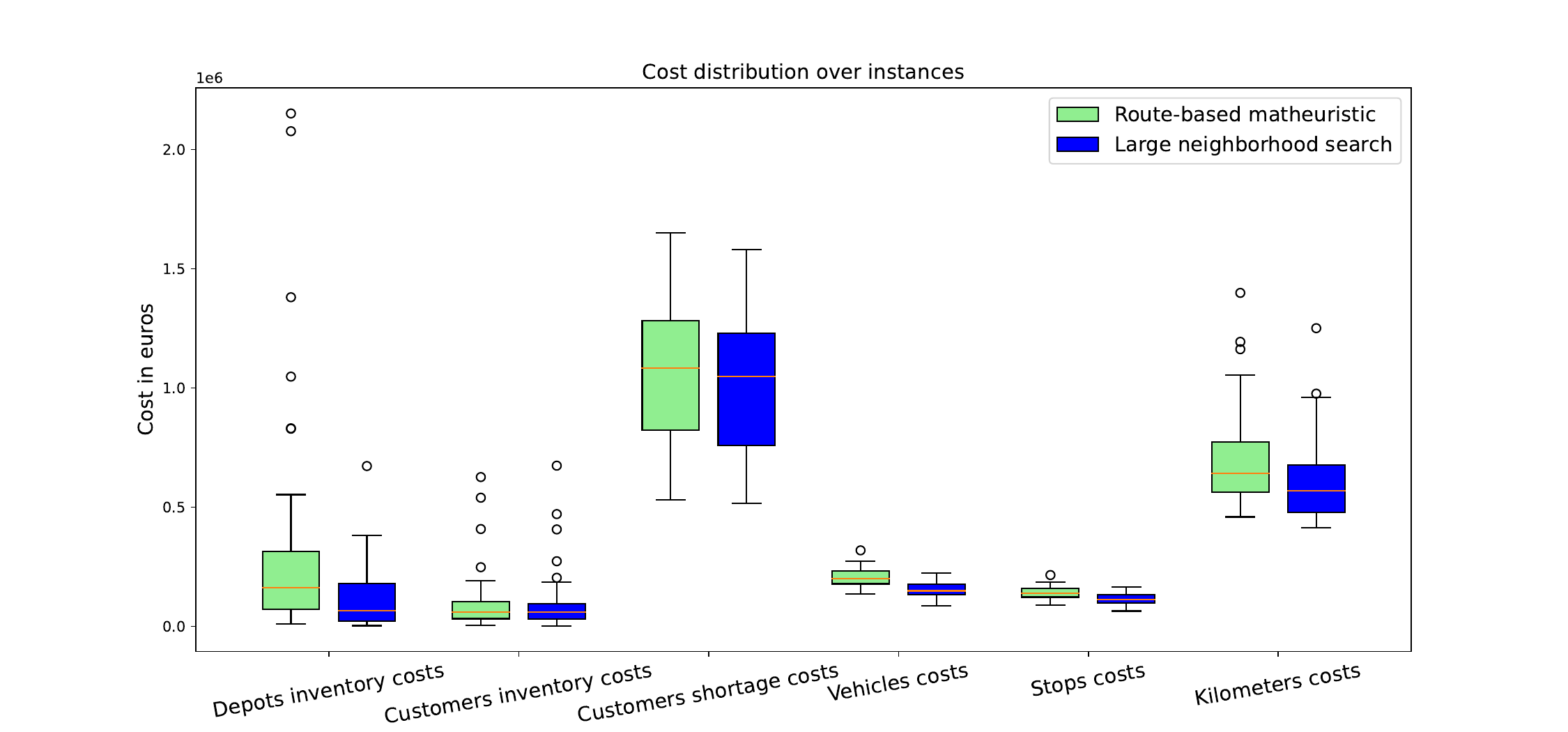}
    \caption{Box plots of the distribution of the solution cost per origin over instances.\\ Green rectangles are related to the 
    route-based matheuristic, blue ones to the LNS.}
    \label{fig:cost_distribution_results}
 \end{figure}

The cost of the IRP solution after our LNS is 
the first natural metric to evaluate our approach.
We can compare it with the cost after the route-based matheuristic defined 
in Section\:\ref{sec:algos_overview}.
On Figure\:\ref{fig:cost_distribution_results}, we show the box plots of the 
cost due to the depots' inventory, the customers' inventory,
the customers' shortage, the vehicles, the stops and the kilometers. The orange 
lines indicate the median over the instances, the ends of the boxes the extreme quartiles ($Q1$ and\:$Q3$),
and the whiskers the range\:$[Q1 - 1.5(Q3-Q1), Q3 + 1.5(Q3-Q1)]$. Outlier points correspond to data outside the whiskers. 
The average total cost over the instances solved by the initialization + local search algorithm is\:$2.88$M€, compared with\:$2.48$M€ after the route-based matheuristic, and\:$2.10$M€ after the LNS.

\begin{result}
  The LNS enables $37\%$ and $18\%$ cost savings compared with the initialization + local search and route-based matheuristic respectively.
 \end{result}

On Figure\:\ref{fig:cost_distribution_results} we notice that almost all the components of the cost are lower: the depots inventory costs, the shortage costs, as well as the routing costs (vehicles costs, stop costs and kilometer costs)
have quartiles corresponding to smaller values after the LNS.

 \begin{result}
  The LNS improves almost every component of the cost.
 \end{result}

\subsection{Extended analysis}\label{subsec:neighborhoods_analysis}

\subsubsection{Time per Operator}\label{subsubsec:gain_time_operators}

\begin{figure}[!h]
  \begin{minipage}{0.55\textwidth}
    \centering
    \includegraphics[width=.99\linewidth]{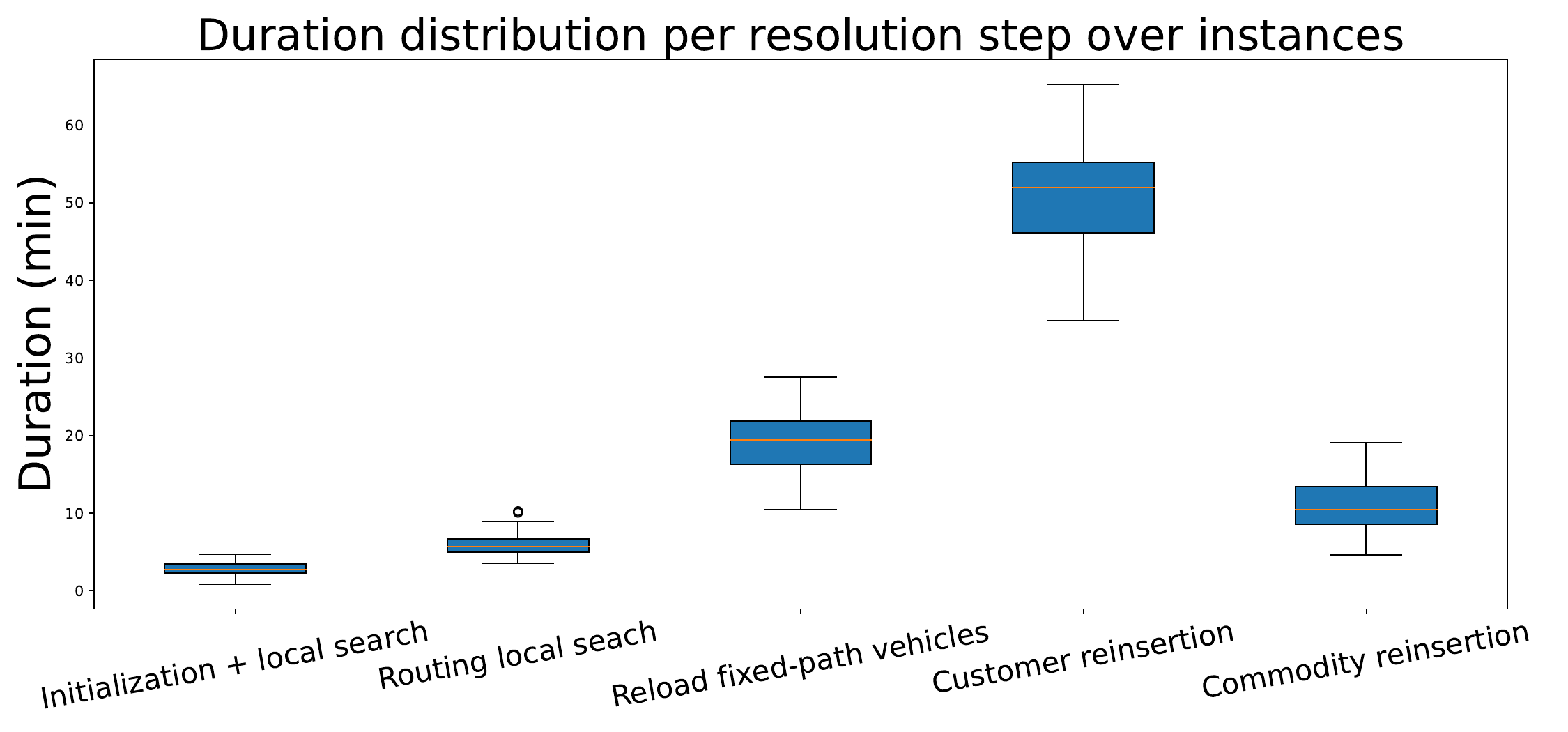}
  \end{minipage}
     \begin{minipage}{0.45\textwidth}
    \centering
    \includegraphics[width=.99\linewidth]{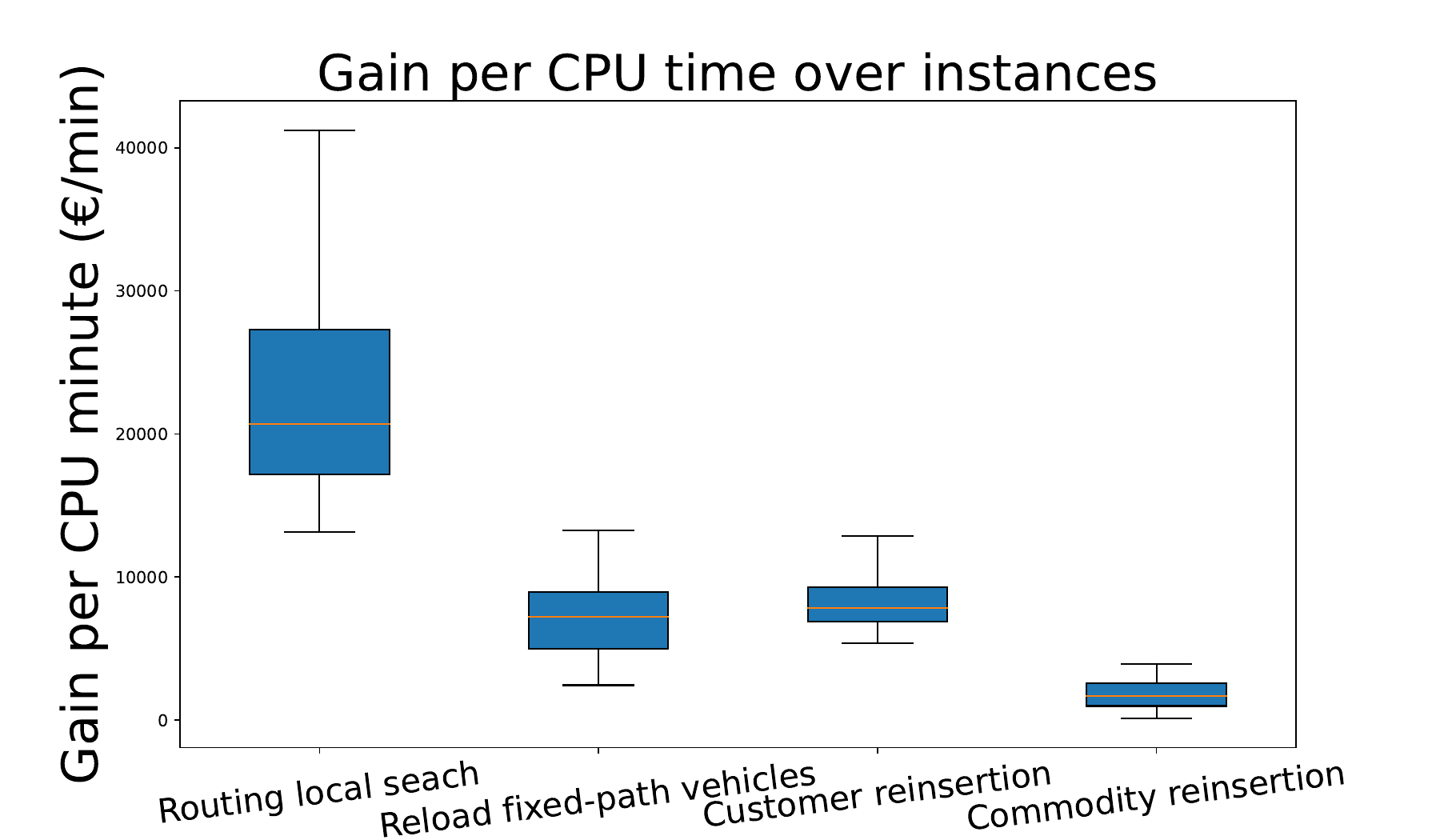}
  \end{minipage}
    \caption{Box plots of the time spent, and of the cost gain per CPU time per operator.}
    \label{fig:duration_gain_per_CPU_distribution_results}
 \end{figure}

We proceed to a more detailed analysis of the solution steps illustrated on Figure\:\ref{fig:algorithms}.
We first compare the total time spent in the initialization + local search algorithm, the routing local search,
the reload fixed-path vehicles large neighborhood, and the customer and commodity reinsertion perturbations.
We emphasize that if we just proceed to one pass over every neighborhood and perturbation without restricting the depots, customers, or commodities considered for the MILP-based operators, we can barely do two outer iterations of the LNS in\:$90$\:minutes.
With the hyperparameters we set, we do\:$2.5$ outer iterations of the LNS on average.
Adaptive approaches \parencite{gendreauHandbookMetaheuristics2010} are thus excluded in our context.
On the left of Figure\:\ref{fig:duration_gain_per_CPU_distribution_results} we show that the initialization + local search algorithm is indeed fast, taking\:$2.89$ minutes on average.
The total time spent in the routing local search amounts to \:$5.94$ minutes on average. 
The reload fixed-path vehicles neighborhood, customer reinsertion and commodity reinsertion perturbations account for
respectively\:$19.2$,\:$51.0$ and\:$11.2$ minutes. About \:$90\%$ of the 
time is spent solving MILPs to modify the solution at a large scale. Therefore, there is not much to gain from improving our implementation further, because the solver itself is out of our control.

\begin{result}
  There is not much performance to gain by further improving our implementation, since most of the time is spent to solve MILPs.
\end{result}

\subsubsection{Cost Gain Over Time}\label{subsubsec:gain_time_ratio}

On the right of Figure\:\ref{fig:duration_gain_per_CPU_distribution_results}, we display the box plots of the cost gain per CPU time (in €/minute)
of our four main LNS components: the routing local search, reload fixed-path vehicles neighborhood, customer reinsertion, and commodity reinsertion. We highlight that the routing local search has the highest gain per CPU time,\:$25.7$k€/minute on average. Its large variance over the instances can be related to the balance 
between routing and inventory costs, which may vary from one instance to another depending on the demand and release profiles. The unit costs themselves do not change between instances.
We thus highlight the routing neighborhoods are efficient and crucial. It is partly due to the structure of our perturbations: a customer or commodity reinsertion 
step creates new routes that are only direct. The reason for this choice is to avoid spending too much time solving greater perturbation MILPs, for they
only involve a portion of the solution. It is therefore useful to combine them with routing neighborhoods that merge or mix routes in different manners to increase route length.

The reload fixed-path vehicles neighborhood, customer and commodity reinsertion perturbations alter both routing and inventory variables, with \:$7.95$k€/minute,\:$8.75$k€/minute and\:$2.77$k€/minute ratios on average. 
Therefore, the descent neighborhoods (routing local search and reload fixed-path vehicles subroutines) decrease the cost on average, but so do the perturbations, even though they are designed to escape from local minima. We show further details with an ablation study in Section\:\ref{subsubsec:ablation_tests}, and tuning in Appendix\:\ref{app:hyperparameters}.
We emphasize that separating the gains per operator is not totally obvious. It is indeed the mix between the axes of the IRP (from routes to commodities and customers) that enables us to substantially reduce the cost in our LNS.

\subsubsection{Ablation study}\label{subsubsec:ablation_tests}

\begin{figure}[htb]
  \includegraphics[width=.9\linewidth]{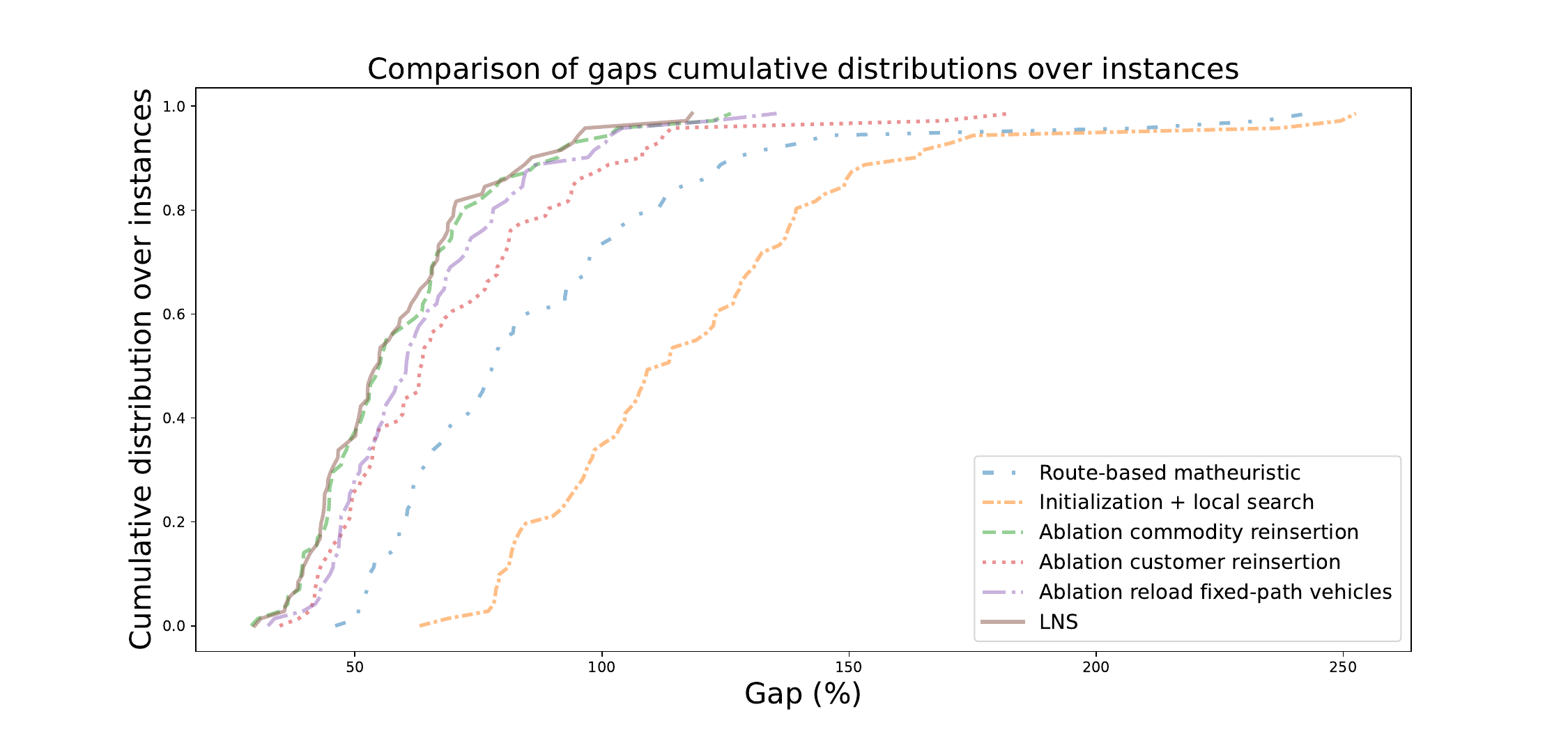}
  \caption{Cumulative distributions of the gap among instances solutions.}
  \label{fig:cumul_gap}
\end{figure}

During the ablation study, we run the LNS with the same hyperparameters, but without the customer reinsertion, commodity reinsertion or reload fixed-path vehicles neighborhood respectively. We compare the cumulative distributions of the gaps on Figure\:\ref{fig:cumul_gap}. We must emphasize our lower bound 
derived in Section\:\ref{min_cost_relaxation} is not tight, and therefore leads to overestimated gaps. It seems wiser to use these gaps as a relative metric to compare the solution processes, and not as an absolute indicator of solution quality. The first remark we can make is that the route-based matheuristic (blue curve) performs better than our initialization + local search algorithm (orange curve), but worse than the LNS, even when one of the LNS components is removed. This confirms the conclusion drawn from the cost analysis in Section\:\ref{subsec:cost_results},
emphasizing the performance of our LNS. Besides, for this short $90$ minutes run, we show that the customer reinsertion perturbation and the reload fixed-path vehicles large neighborhood are crucial for the LNS performance. Indeed, removing the customer reinsertion subroutine (red curve) entails a $12\%$ gap increase on average. Removing the reload fixed-path vehicles leads to $6\%$ average increase (purple curve). Those two subroutines are complementary 
% seem to be relatively ``orthogonal'' 
in the sense that they modify distinct parts of the solution. The commodity reinsertion perturbation seems to be less crucial for the performance in the short $90$ minutes run. Indeed, removing it from the LNS entails only a $2\%$ average gap increase. It may be due to the fact that changing the commodity flow axis is also at the core of the reload fixed-path vehicles subroutine. Nonetheless, the commodity reinsertion perturbation potentially creates new direct routes. We show additional experiments in Appendix\:\ref{app:hyperparameters} with longer runs to emphasize its utility.

\begin{result}
  All the neighborhoods are useful. The commodity reinsertion perturbation only becomes fully efficient for longer runs.
\end{result}

\subsubsection{Results with Longer Routes and Larger CPU Time}\label{subsubsec:res_longer_routes}

As stated above, limiting the routes to\:$S_{\max} = 3$ stops maximum is restrictive, though this is an important constraint for Renault. To emphasize the benefits induced by longer routes, we solve the instances with\:$S_{\max} \in \{3,10\}$, and a time limit of\:$300$ minutes. Note that both the initialization + local search and route-based matheuristic end before the time limit. The rest of the hyperparameters are unchanged. We summarize the average cost impact with Table\:\ref{tab:smaxvary}. 
First, increasing the LNS running time from $90$ minutes to $300$ minutes with $S_{\max} = 3$ reduces the average gap by an additional seven percent. More experiments on cost evolution with LNS iterations are analyzed in Appendix\:\ref{app:cost_evol}. There we also show that we can reach lower costs with $S_{\max} = 10$ with each algorithm, observing a large gain after the initialization + local search algorithm. Numerical experiments also demonstrate that our LNS can handle the combinatorics induced by longer routes.

\begin{result}
  Increasing the running time enables to find better solutions. The LNS has not converged after $90$ minutes.
\end{result}

\begin{table*}[!h]
  \centering
  \caption{Comparison of the average costs over instances between\:$S_{\max} = 3, 10$ after $300$ minutes.} 
  \label{tab:smaxvary} 
  \setlength{\tabcolsep} {0.3cm} 
  \scalebox{0.92} 
  {   
    \begin{tabular}{ccccc}
  \toprule
                              & \begin{tabular}{@{}c@{}}\textbf{Initialization +} \\ \textbf{local search}\end{tabular} & \begin{tabular}{@{}c@{}}\textbf{Route-based} \\ \textbf{matheuristic}\end{tabular}  & \textbf{LNS} &\textbf{Lower bound}\\
  \midrule
  Average cost $S_{\max} = 3$ & $2.88$M€ &  $2.48$M€ & $2.01$M€ & $1.35$M€\\
  Average cost $S_{\max} = 10$& $2.41$M€ &  $2.22$M€ & $1.81$M€ & $1.35$M€\\
  \bottomrule
  \end{tabular} %
  }
\end{table*}

\begin{result}
  Increasing the maximum length of a route\:$S_{\max}$ enables additional savings when running time can be increased.
\end{result}

\section{Conclusions}\label{conclusion}

In this study motivated by an industrial partnership with Renault, we consider European-scale continuous-time multi-attribute IRP
instances. This inherently hard problem has to be solved in a limited time of\:$90$ minutes every day.
To do so, we design a large neighborhood search based on the generalization of TSP and SDVRP neighborhoods Section\:\ref{sec:routing_LS}, 
a large neighborhood Section\:\ref{reload} inspired by recent matheuristics designed for the IRP \parencite{bertazziMatheuristicAlgorithmMultidepot2019,archettiMatheuristicMultivehicleInventory2017}, 
and two new perturbations Section\:\ref{customer_reinsertion} and Section\:\ref{commodity_reinsertion} based on MILPs. 
We also derive an initialization + local search algorithm that relies on flows Section\:\ref{min_cost_relaxation}. It allows us to quickly initialize our 
IRP instances with non-trivial solutions, and to derive a lower bound. To the best of our knowledge, this lower bound is unknown in the IRP literature.
We extract and process a dataset of\:$71$ European-scale multi-attribute IRP instances that we make available publicly.
Numerical experiments in Section\:\ref{sec:numerical_results} show the results of our LNS.
We highlight that it outperforms the route-based matheuristic, and that each component of the LNS brings useful contributions.

We also emphasize some limits and perspectives. The choice of the neighborhoods is hard to define a priori, and we do not have time to browse the space of neighborhoods repeatedly. Parallel computing and multi-threading could be some perspectives to consider, although even small changes imply inventory dynamics over the whole time horizon, leading to overlaps between neighborhoods. Cache use would thus be a challenge. Some techniques of machine learning for operations research could also be considered in this direction, for example in a reinforcement learning \parencite{wulearninglargeneighborhoodsearchpolicy} or structured learning \parencite{parmentierLearningApproximateIndustrial2022} paradigm. Besides, we suffer from the poor quality of our lower bound to compute gaps. Deriving a better relaxation for this multi-attribute IRP is a challenge that has not been addressed in the literature, to the best of our knowledge. We could add valid 
inequalities leveraging the research on exact solutions of smaller problems \parencite{manousakisImprovedBranchandcutInventory2021,desaulniersBranchPriceandCutAlgorithmInventoryRouting2015} to improve the quality of our neighborhoods. A version of the LNS we designed is currently used in production at Renault, with additional statistical treatments due to the noise observed in release and demand forecasts. Our simulations enable us to estimate a reduction of thousands of tons of $\text{CO}_2$ and millions of euros per year. 
%We hope to confirm it with real savings within the next months.
Dealing with the stochastic IRP \parencite{nolzStochasticInventoryRouting2014,coelhoDynamicStochasticInventoryrouting2012} is at the core of our future work.

% Acknowledgments here
\paragraph*{Acknowledgments}
We are grateful to the Renault supply-chain and IT teams for the partnership we have set and the industrial motivation they provide us with, especially to Alain Nguyen, Thaddeus Leonard, Nicusor-Eugen Plescan, Christian Serrano, Ludovic Doudard and Aimé-Frédéric Rosenzweig. We also would like to thank Vincent Leclère for his advice on the form of the article.

\printbibliography

\begin{appendix}
  
\section{Notations}\label{app:notations}
In Table\:\ref{data_table}, we recap the main concepts and 
notations for the dimensions, the initial inventory, the demand and release, 
the free inventory capacities, the locations graph, the commodities and vehicles lengths, and the unit costs.

\begin{table*} [!t]
\centering
\caption{Notations} \label{data_table}
\setlength{\tabcolsep} {0.3cm} 
\scalebox{0.95} 
{ 
\begin{tabular} {ll} 
\toprule
\textbf{Name} & \textbf{Description} \\
\midrule
$t$                     & Day                                                                         \\
$T$                     & Horizon                                                                      \\
$d$                     & Depot (facility)                                                              \\
$D$                     & Set of depots     \\
% $w$ & plateforme \\
% $W$ & Ensemble des plateformes \\
$c$                     & Customer                                                                  \\
$C$                     & Set of customers                                                    \\
% $\tilde F$ & Ensemble des fournisseurs livrés en direct \\
% $F^p$ & Ensemble des fournisseurs livrés par $w$ \\
$m$                     & Commodity                                                             \\
$M$                     & Set of commodities \\ 
$M_c$                   & Set of commodities used by customer $c$                      \\
$M_d$                   & Set of commodities used by depot $d$ \\
% $E_{u,f}$ & Ensemble des emballages dédies à $(u,f)$ \\
\midrule
$\bfI_{md}^0$               & Initial inventory of $m$ at $d$ \\
$\bfI_{mc}^0$               & Initial inventory of $m$ at $c$ \\
$b_{mdt}^+$             & Number of commodity $m$ released on day $t$ by $d$          \\
$b_{mct}^-$             & Number of commodity $m$ demanded on day $t$ by $c$                         \\
% $b_{efj}^+$ & Nombre de piles de $e$ livrées jour $j$ par $f$\\
% ${B_{efj}^+}$ & Nombre maximum de piles/palettes de $e$ à livrer à $f$ entre $1$ et $j$ \\
% $\underline{B_{efj}^+}$ & Nombre minimum de piles/palettes de $e$ à livrer à $f$ entre $1$ et $j$ \\
${\kappa_{mdt}}$             &  Free inventory capacity of $m$ at $d$ on the evening of day $t$                     \\
% ${\kappa_{ewjj}}$ & Stock maximum de piles de $e$ en $w$ le soir de $j$ \\
% $\underline{\kappa_{efj}}$ & Stock minimum contractuel de $e$ en $f$ le soir de $j$ \\
${\kappa_{mct}}$             &  Free inventory capacity of $m$ at $c$ on the evening of day $t$                 \\
\midrule
% $L$ & Métrage linéaire d'un camion \\
$\mathcal{V}$                     & Set $D \cup C$ of the vertices of the locations graph                          \\
$v$                         & Node of $\mathcal{V}$                                                               \\
$\mathcal{A}$                     & Arcs of the locations graph \\%, $(d \to c)$ or $(c_1 \to c_2)$ for $d \in D$, $c, c_1, c_2 \in C, c_1 \neq c_2$                              \\
$a$                         & Arc                                                                        \\
$\Delta_a$                   & Distance in kilometers corresponding to arc $a$                            \\
$\tau_a$                     & Duration in hours corresponding to arc $a$                                  \\
% $c_a$ & Coût à opérer l'arête $a$ avec un camion \\
% $t_a$ & Durée nécessaire pour parcourir l’arête $a$ en camion \\
$P$                     & Path in the locations graph                                                       \\
$S_{max}$               & Maximum number of stops in a route                    \\
$\tau_{\max}$                  & Number of transport hours per day \\
% $\calR^{\mathrm{r}}$ & Ensemble des routes régulières \\
% $\calR^{\mathrm{s}}$ & Ensemble des routes spot \\
% $\calR_{u,f}^{\mathrm{d}}$ & Ensemble des routes dédiées à $u,f$. \\
$\ell_m$                & Length of a commodity $m$                                           \\
$L$                     & Length of a vehicle (homogeneous)                                                                         \\
% $M^{\mathrm{r}}(j,P) $& Nombre maximum de routes régulières partant le jour $j$ et suivant le chemin $P$ \\
% $c^{\mathrm{t}}(j,P)$& Coût négocié d'un route régulière partant le jour $j$ et suivant le chemin $P$ \\
% $C$ & Cluster géographique \\
% $\calC_u$& Partition des fournisseurs servis par $u$ en clusters géographiques \\
\midrule
$c_{md}^{\texttt{exc}}$   & Unit excess inventory storage cost of $m$ at $d$ per night                 \\
% $c_{ew}^{\mathrm{s}}$ & Coût unitaire d'un stock excédentaire de $e$ en $w$ par nuit \\
$c_{mc}^{\texttt{exc}}$   & Unit excess inventory storage cost of $m$ at $c$ per night          \\
% $c_{ef}^{\mathrm{nrc}}$ & Coût journalier de non respect du contrat avec $f$ pour $e$ \\
$c_{mc}^{\texttt{short}}$ & Unit shortage cost of $m$ at $c$ per day\\
$c^{\texttt{km}}$                & Per kilometer cost of the routes                                                \\ %spot du cluster géographique $C$ \\
$c^{\texttt{veh}}$      & Unit cost for using a vehicle                                                             \\
$c^{\texttt{stop}}$     & Unit cost for making a stop                                                  \\
\bottomrule
\end{tabular} %
} 
\end{table*}

\section{Reinsertion cost}\label{app:reinsertion_cost}

In Table\:\ref{tab:arcs_vehicle_graph_customer_reinsertion}, we introduce\:$c^r_{sc}$, the cost induced by the insertion of customer\:$c$ in route\:$r$ at position\:$s$. We define it explicitly here.
Let\:$r$ be a route, with day of departure $t^r$, path $P^r = (v_0^r, v_1^r, \ldots, v_k^r)$ and quantities delivered ${\mathbf{q}^r= (q^r_{ms})_{m \in M, s \in [|P^r|-1]} \in (\mathbb{Z}^+)^{|M|\times(|P^r|-1)}}$. Let $c$ be a customer which does not belong to $P^r$ and $s$ its insertion position. We denote by\:$\tilde{r}$ the route with corresponding insertion, $P^{\tilde{r}}$ the associated path, and\:$\tilde{I}$ the inventory variables with updated delivery times due to delays induced by the new stop (see Section 2.3. of the paper). The cost\:$c^r_{sc}$ is induced by the insertion of customer\:$c$ at position\:$s$ in the stops of route\:$r$. It involves routing and inventory considerations at the other customers delivered by route\:$r$, because of the delays. It is computed as:
  \begin{align*}
  & c^r_{sc} = c^{\texttt{stop}} + c^{\texttt{km}} \bigg[\Delta_{(c_{s-1}, c)} + \Delta_{(c, c_{s+1})} - \Delta_{(c_{s-1}, c_{s+1})} \bigg] + 
  \\ & \sum_{\substack{c' \in (v_1^r, \ldots, v_k^r), \\t \geq t^r, \\m, \sum_{s \in [|P^r|-1]} q^r_{sm} > 0 } } \bigg[c^{\texttt{exc}}_{mc'}\left(\tilde{I}_{mc't} - \kappa_{mc't}\right)^+ - c^{\texttt{exc}}_{mc'}\left(I_{mc't} - \kappa_{mc't}\right)^+ \bigg] +
  \\ & \sum_{\substack{c' \in (v_1^r, \ldots, v_k^r), \\t \geq t^r, \\m, \sum_{s \in [|P^r|-1]} q^r_{sm} > 0 } } \bigg[c^{\texttt{short}}_{mc'} (b^-_{mc't} -\tilde{I}_{mc'(t-1)})^+ - c^{\texttt{short}}_{mc'} (b^-_{mc't} -I_{mc'(t-1)})^+ \bigg].
  \end{align*} 
  It has one term of routing cost, one of excess inventory cost, and one of shortage cost. The two last terms are reduced to the customers visited by the route\:$r$, and browsing only days and commodities impacted by route\:$r$. This local implementation is important for code performance. Besides, the value of\:$c^r_{sc}$ does not depend on the commodity flows at the inserted customer, contrary to the inventory cost at the starting depot of the routes, and the inventory cost at the inserted customer. It can therefore be pre-computed before solving the customer reinsertion MILP.

\section{Instance details.}\label{app:instance}
\begin{figure}[!htb]
    \includegraphics[width=.9\linewidth]{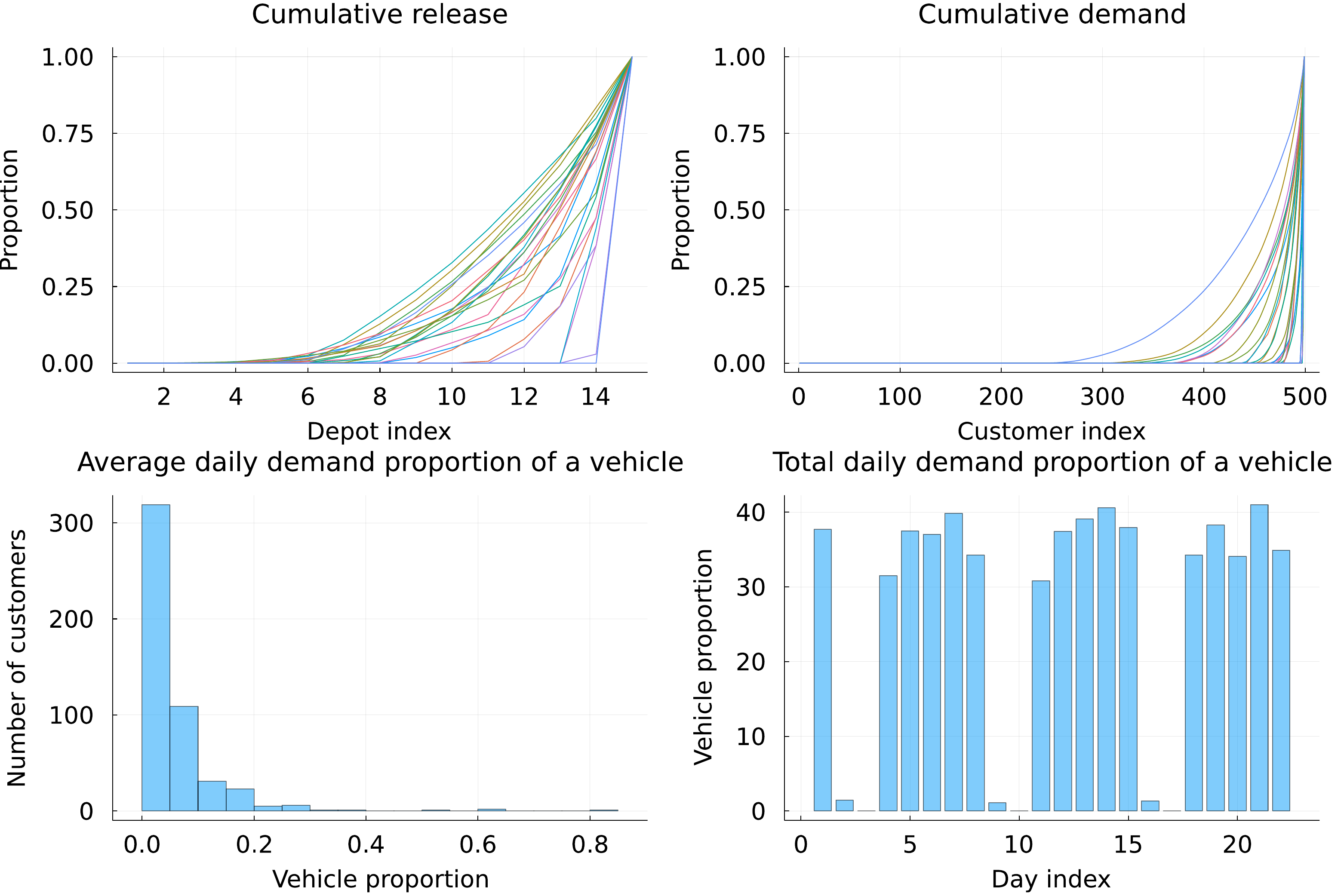}
    \caption{Details on the release and demand features of an instance.}
    \label{fig:instance}
  \end{figure}

  The complexity of the industrial problem for Renault comes from the tight coupling between commodities, between nodes (depots and customers), between inventory and routing (increased by the continuous-time aspect), as well as the size of the instances. We provide some details about the instances' structure here to back our discussion. In Figure\:\ref{fig:instance}, we focus on one particular instance with $15$ depots and $499$ customers. The code to generate the plots, as well as the instances are publicly available. On the top-left plot, we show the cumulative normalized total release distribution over depots per commodity. For each commodity (each line in the plot), we sort depots per total release. Therefore, the depot index axis is different for each commodity. Our aim here is to show how the release of each commodity is spread over depots. We see that for a few commodities, only two or three depots release almost the whole proportion. But for most of the commodities, the release is shared by about eight depots. We recall the distribution of the number of depots per instance is visible in Figure\:\ref{fig:dimensions_instances}. Similarly, on the top-right plot, we show the cumulative normalized total demand distribution over customers per commodity. We see that no commodity is shared by the $499$ customers. One is shared by more than $200$ customers, and most of them are shared by $50$ to $100$ customers. On the bottom-left plot, we show the histogram over customers of the fraction of a vehicle represented by the average daily demand. On each day and for each customer, the commodities are scaled by their lengths to compute the total demand length of the day, and the daily average is computed. In this way, we see that a huge proportion of the customers have a demand that fits in less than $10\%$ of a vehicle. Only one customer has a demand that requires $80\%$ of a vehicle. This analysis has to be combined with the last plot of Figure\:\ref{fig:instance}. On the bottom-right plot, we show the proportion of a vehicle length represented by the total daily demand length. Apart from the weekend, we observe that the total demand length represents between $30$ and $40$ times the length of a vehicle. The bin packing, release profiles and continuous-time aspects lead to many more vehicles per day in practice, as illustrated in Figure\:\ref{fig:solution}. We thus show it is crucial to create routes with several stops to find good solutions. 

  \begin{result}
    Instances are complex because of the bindings between depots, customers, days and commodities. 
  \end{result}

\section{Solution details.}\label{app:solution}
\begin{figure}[!htb]
    \includegraphics[width=.9\linewidth]{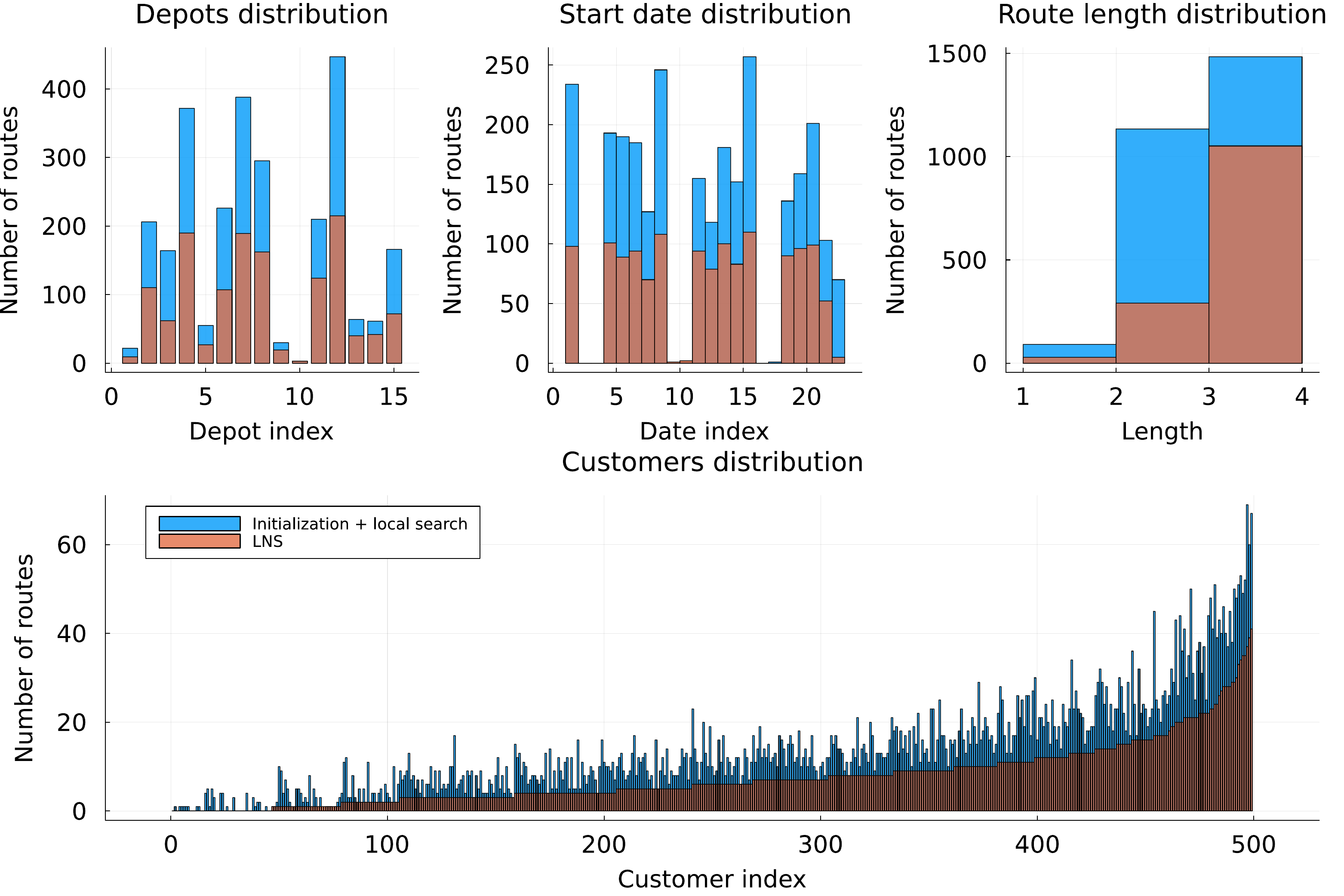}
    \caption{Details on LNS and initialization + local search solutions of an instance.}
    \label{fig:solution}
  \end{figure}

  We focus here on the best solution found by the LNS after a $90$ minutes run, and compare it to the initialization + local search solution, for the instance highlighted in the previous appendix. As for the instance analysis, the code to generate the plot is publicly available, as well as the solutions found. On Figure\:\ref{fig:solution}, we display depot, start date, number of stops and customer distributions over the routes of the two solutions. Blue color is related to the initialization + local search, red color to the LNS. The first remark we can make is that there are much fewer routes in the solution of the LNS than in the solution of the initialization + local search, which explains the decrease of the total routing cost after the LNS. On the top-left of Figure\:\ref{fig:solution}, we see that both the solution of the initialization + local search and the best one found by the LNS have routes starting from every depot. The distribution structures are similar, although the numbers of routes are distinct. This can be linked to the release and demand structures. On the top-middle plot, we show that except for the weekends, the departure dates are spread over days. The LNS solution has less starting date variations than the initialization + local seach solution. We expect it to be linked to a better management of the time axis, with the challenging couplings between inventory dynamics and continuous-time routes. We also observe the effect of the end of the time axis, with a decreasing number of departures when approaching the horizon. On the top-right plot of Figure\:\ref{fig:solution}, we see that only a few routes in both solutions are direct, meaning they visit an only one customer. We also observe that the LNS leads to a much larger proportion of routes with three stops than the initialization + local search. Recall the remark made on Figure\:\ref{fig:instance} in the previous appendix: the demand distribution is spread over customers, so there is a need to group them in long routes. Last, the bottom plot on Figure\:\ref{fig:solution} shows the distributions of the number of visits per customer in both the solution of the LNS and of the initialization + local search. We see that visits are widely spread over customers. The LNS seems to cancel some visits to customers with small numbers of visits after the initialization + local search, those visits may be related to very expensive routes. 

  \begin{result}
    The strong couplings between depots, customers, days and commodities is also found in the solutions structure. 
  \end{result}

\section{LNS hyperparameters tuning.}\label{app:hyperparameters}
  
  In the LNS, we jointly tune two hyperparameters: the number of customer reinsertion steps, and the number of commodity reinsertion steps (both for one iteration of the LNS). Their values can range from zero to the number of customers or number of commodities in the instance solved. We could also tune the number of reload fixed-path vehicles steps, as well as the number of routing local search steps per outer LNS iteration. The space of configurations to consider would become too huge to be explored with reasonable computing resources. We thus focus here on the perturbations hyperparameters. We display results in Table\:\ref{tab:gridsearch}: the mean, maximum and minimum gaps over the $71$ instances after $90$ and $300$ minutes runs when varying the hyperparameters, as well as the average number of outer LNS steps done after $90$ and $300$ minutes. Table\:\ref{tab:gridsearch} also displays the statistics of the initialization + local search algorithm as a reference. Overall, we see that for $90$ minutes runs, increasing the number of customer reinsertion steps is crucial to achieve better performance. As shown in the ablation study in Section\:\ref{subsubsec:ablation_tests}, the commodity reinsertion perturbation is less important. Yet, the best results are found when doing ten commodity reinsertion steps, and browsing the whole set of customers for the customer reinsertion subroutine per LNS outer steps. It leads to a $58\%$ average gap in $90$ minutes. For longer runs of $300$ minutes, the commodity reinsertion perturbation is more important. Indeed, the worst hyperparameters configurations found are those with zero commodity reinsertion step. We thus highlight this perturbation is useful in longer runs, when escaping from local minima becomes crucial. 

  \begin{result}
    Our hyperparameters are tuned for short runs. Longer runs highlight the commodity reinsertion perturbation is relevant to escape from local minima.
  \end{result}

  \begin{table*}[!ht]
    \centering
    \caption{Gap results when tuning the numbers of customer and commodity reinsertion steps per outer LNS step on the $71$ instances.}
    \setlength{\tabcolsep} {0.2cm} 
  \scalebox{0.6} 
    { 
    \begin{tabular}{ccccccc}
    \toprule
        \textbf{Algorithm} & \begin{tabular}{@{}c@{}}\textbf{Commodity} \\ \textbf{reinsertion steps}\end{tabular} & \begin{tabular}{@{}c@{}}\textbf{Customer} \\ \textbf{reinsertion steps}\end{tabular} & \begin{tabular}{@{}c@{}}\textbf{Mean} \\ \textbf{gap}\end{tabular} & \begin{tabular}{@{}c@{}}\textbf{Min} \\ \textbf{gap}\end{tabular} & \begin{tabular}{@{}c@{}}\textbf{Max} \\ \textbf{gap}\end{tabular} & \begin{tabular}{@{}c@{}}\textbf{Mean number} \\ \textbf{of LNS steps}\end{tabular}\\ 
        \midrule
        \begin{tabular}{@{}c@{}}\textbf{Initialization +} \\ \textbf{local search}\end{tabular} & 0 & 0 & 119 & 63 & 252 & 0 \\ 
        \midrule
        &0 & 0 & 70 & 35 & 160 & 12,8 \\  
        &5 & 0 & 69 & 35 & 164 & 8,18 \\  
        &10 & 0 & 70 & 35 & 183 & 6,03 \\  
        &15 & 0 & 70 & 36 & 179 & 4,48 \\  
        &20 & 0 & 71 & 36 & 197 & 4,07 \\  
        &25 & 0 & 71 & 36 & 183 & 3,07 \\  
        &max & 0 & 72 & 37 & 190 & 2,72 \\  
        &0 & 50 & 63 & 32 & 139 & 10,9 \\  
        &5 & 50 & 64 & 34 & 143 & 7 \\  
        &10 & 50 & 63 & 34 & 135 & 5,44 \\  
        &15 & 50 & 65 & 34 & 154 & 4,32 \\  
        &20 & 50 & 64 & 35 & 145 & 3,7 \\  
        &25 & 50 & 66 & 35 & 154 & 3,06 \\  
        &max & 50 & 65 & 36 & 151 & 2,87 \\  
        &0 & 100 & 62 & 31 & 140 & 9,16 \\  
        &5 & 100 & 61 & 32 & 132 & 6,44 \\  
        &10 & 100 & 62 & 33 & 133 & 4,83 \\  
        &15 & 100 & 62 & 32 & 129 & 4,04 \\  
        &20 & 100 & 64 & 34 & 136 & 3,42 \\  
        &25 & 100 & 65 & 34 & 161 & 2,75 \\  
        &max & 100 & 65 & 35 & 132 & 2,55 \\  
        &0 & 200 & 60 & 31 & 131 & 7,37 \\  
        &5 & 200 & 61 & 32 & 127 & 5,44 \\  
        \begin{tabular}{@{}c@{}}\textbf{LNS} \\ \textbf{$90$ minutes}\end{tabular}&10 & 200 & 61 & 32 & 127 & 4,17 \\  
        &15 & 200 & 61 & 32 & 121 & 3,64 \\  
        &20 & 200 & 62 & 33 & 146 & 2,9 \\  
        &25 & 200 & 61 & 33 & 141 & 2,63 \\  
        &max & 200 & 62 & 34 & 126 & 2,32 \\  
         & 0 & 300 & 60 & 30 & 125 & 5,95 \\  
        &5 & 300 & 59 & 31 & 127 & 4,8 \\  
        &10 & 300 & 60 & 31 & 120 & 3,73 \\  
        &15 & 300 & 61 & 32 & 125 & 2,92 \\  
        &20 & 300 & 60 & 32 & 120 & 2,9 \\  
        &25 & 300 & 62 & 33 & 124 & 2,24 \\  
        &max & 300 & 63 & 33 & 125 & 1,85 \\  
        &0 & 400 & 59 & 31 & 128 & 4,96 \\  
        &5 & 400 & 60 & 31 & 127 & 3,76 \\  
        &10 & 400 & 59 & 31 & 120 & 3,11 \\  
        &15 & 400 & 60 & 33 & 131 & 2,55 \\  
        &20 & 400 & 61 & 31 & 123 & 2,39 \\  
        &25 & 400 & 59 & 33 & 131 & 1,91 \\  
        &max & 400 & 62 & 34 & 121 & 1,86 \\  
        &0 & max & 60 & 29 & 126 & 3,55 \\  
        &5 & max & 59 & 31 & 118 & 2,97 \\  
        &\textbf{10} & \textbf{max} & \textbf{58} & \textbf{30} & \textbf{118} & \textbf{2,54} \\  
        &15 & max & 59 & 30 & 122 & 2,27 \\  
        &20 & max & 59 & 30 & 123 & 1,9 \\  
        &25 & max & 59 & 32 & 120 & 1,77 \\  
        &max & max & 59 & 31 & 131 & 1,63 \\  
        \midrule 
        & 0 & 200 & 55 & 28 & 116 & 31.1 \\  
        & 10 & 200 & 52 & 27 & 107 & 19.4 \\  
        & max & 200 & 54 & 28 & 117 & 10.1 \\  
        &0 & 400 & 55 & 27 & 127 & 24.2 \\  
        \begin{tabular}{@{}c@{}}\textbf{LNS} \\ \textbf{$300$ minutes}\end{tabular}&\textbf{10} & \textbf{400} & \textbf{51}  & \textbf{26} & \textbf{104} & \textbf{15.9} \\  
        &max & 400 & 53 & 27 & 113 & 9.58 \\  
        & 0 & max & 55 & 27 & 117 & 21.7 \\  
        &10 & max & 52 & 26 & 110 & 12.2 \\  
        &max & max & 51 & 26 & 115 & 8.15 \\  

        \bottomrule
    \end{tabular}\label{tab:gridsearch} 
    }
\end{table*}
  
\section{LNS cost evolution.}\label{app:cost_evol}
\begin{figure}[!htb]
    \includegraphics[width=.9\linewidth]{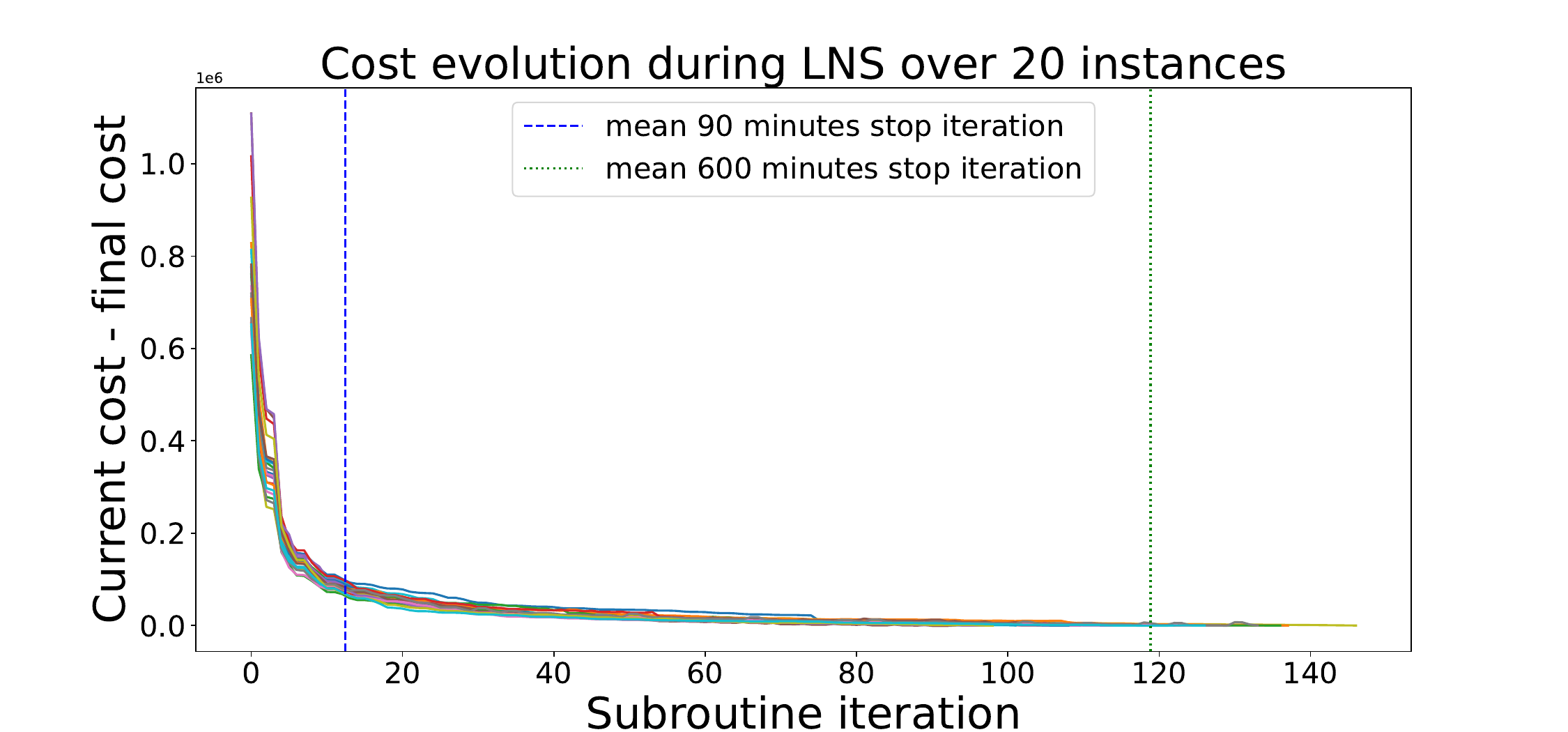}
    \caption{Evolution of the LNS cost over $20$ instances during $600$ minutes runs.}
    \label{fig:cost_evol}
  \end{figure}

  To extend the analysis made on the longer $300$ minutes runs in Section\:\ref{subsubsec:res_longer_routes}, showing that the LNS has not converged in $90$ minutes, we track cost with LNS iterations in long $600$ minutes runs over $20$ instances. We show in Figure\:\ref{fig:cost_evol} the difference between the current and final costs over subroutine iterations in the LNS. We denote by subroutine iteration the pass in one of the routing local search, reload fixed-path vehicles, customer reinsertion and commodity reinsertion subroutines (possibly making several inner iterations in each of them). Since each subroutine has a distinct behavior on each instance, the same time limit does not lead to the same number of subroutine iterations. Overall, we observe a fast cost reduction phase before $90$ minutes, and a slower phase after. On average, we reduce by $4\%$ the cost after $600$ minutes compared to the $90$ minutes runs.

  \begin{result}
    Although the LNS has not converged in $90$ minutes, cost reduction is much slower on average after this point.
  \end{result}

\end{appendix}

\end{document}